\documentclass[10pt,fleqn]{article}

\newif\ifpdf
\ifx\pdfoutput\undefined
    \pdffalse       
\else
    \pdfoutput=1    
    \pdftrue
\fi

\usepackage{amsmath,amssymb,amsthm}

\ifpdf
    \usepackage[pdftex]{graphicx}
    \usepackage[backref, colorlinks=true, pdfstartview=FitV, linkcolor=blue,
        citecolor=blue, urlcolor=blue]{hyperref}
\else
    \usepackage{graphicx}
    \usepackage{hyperref}
\fi
\usepackage[numeric]{amsrefs} 

    \newcommand{\MathReview}[1]{~\href{http://www.ams.org/mathscinet-getitem?mr=#1}{\mbox{\bf MR~#1}}}

    \newcommand{\NN}{\ensuremath{\mathbb N}}
    \newcommand{\ZZ}{\ensuremath{\mathbb Z}}

    \newcommand{\EE}{\ensuremath{\mathbb E}}
    \newcommand{\PP}{\ensuremath{\mathbb P}}

    \newcommand{\Prob}[1]{\ensuremath{\PP\left[ #1 \right]}}
    \newcommand{\Expect}[1]{\ensuremath{\EE\left[ #1 \right]}}

    \newcommand{\cE}{\ensuremath{{\cal E}}}
    \newcommand{\cF}{\ensuremath{{\cal F}}}
    \newcommand{\cV}{\ensuremath{{\cal V}}}
    
    \newcommand{\cS}{\ensuremath{{\cal S}}}

    \newcommand{\cT}{\ensuremath{{\cal T}}}
    \newcommand{\cP}{\ensuremath{{\cal P}}}
    \newcommand{\cQ}{\ensuremath{{\cal Q}}}
    \newcommand{\cR}{\ensuremath{{\cal R}}}

    \newcommand{\floor}[1]{\mbox{$\left\lfloor #1 \right\rfloor$}}

    \newcommand{\field}[1]{{\mathbb F}_{ #1 }}

    \newcommand{\citebarfnitem}[1]{\ref{lem:barfn}({\em#1}\,)}

    \DeclareMathOperator{\ord}{ord}
    
    \DeclareMathOperator{\lcm}{lcm}

    \newtheorem{thm}{Theorem}[section]
    \newtheorem{lem}[thm]{Lemma}
    
    \newtheorem{cor}[thm]{Corollary}
    
    \newtheorem{cnj}[thm]{Conjecture}
    
    \newtheorem{prop}[thm]{Proposition}

\title{Reciprocals of Binary Power Series}
\author{
    Joshua N. Cooper\thanks{Supported by NSF grant DMS-0303272}\\
    {\small Courant Institute of Mathematics, New York University, New York, NY}\\
    {\small \tt{jcooper@cims.nyu.edu}}\\
    Dennis Eichhorn\\
    {\small California State University, East Bay, Hayward, CA}\\
    {\small \tt{eichhorn@mcs.csueastbay.edu}}\\
    Kevin O'Bryant\thanks{Supported by NSF grant DMS-0202460.}\\
    {\small University of California, San Diego, San Diego, CA}\\
    {\small City University of New York, College of Staten Island, New York, NY}\\
    {\small \tt{kevin@member.ams.org}}}
\date{\today}

\begin{document}
    \ifpdf
       \DeclareGraphicsExtensions{.pdf,.jpg,.mps,.png}
    \fi
\maketitle \sloppypar

\begin{abstract}
If $A$ is a set of nonnegative integers containing 0, then there is a unique nonempty set $B$ of nonnegative
integers such that every positive integer can be written in the form $a+b$, where $a\in A$ and $b\in B$, in an
even number of ways. We compute the natural density of $B$ for several specific sets $A$, including the
Prouhet-Thue-Morse sequence, $\{0\}\cup\{2^n\colon n\in\NN\}$, and random sets, and we also study the
distribution of densities of $B$ for finite sets $A$. This problem is motivated by Euler's observation that if
$A$ is the set of $n$ that have an odd number of partitions, then $B$ is the set of pentagonal numbers
$\{n(3n+1)/2 \colon n\in\ZZ\}$. We also elaborate the connection between this problem and the theory of de
Bruijn sequences and linear shift registers.
\end{abstract}

    \markright{Reciprocals of Binary Power Series}
    \pagestyle{myheadings}
    \thispagestyle{empty}

\section{Introduction}\label{sec:Introduction}

There is a unique set $B$ of nonnegative integers with the property that each positive integer can be written in
the form $s^2+b$ ($s\in \NN:=\{0,1,2,\dots\}, b\in B$) in an even number of ways. Specifically,
    \[
    B= \{ 0, 1, 2, 3, 5, 7, 8, 9, 13, 17, 18, 23, 27, 29, 31, 32, 35, \ldots \}.
    \]
Are the even numbers in $B$ exactly those of the form $2k^2$? Does $B$ have positive density?

Before addressing these two questions, we restate and motivate the problem in greater generality. Given any sets
$A,B\subseteq \NN$, the asymmetric additive representation function is defined by
    \[
    R(n) := \#\{(a,b) \colon n=a+b, a\in A, b\in B\};
    \]
equivalently, we could define $R$ by noting that
    \[
    \bigg( \sum_{a\in A} q^a \bigg) \, \bigg(\sum_{b\in B} q^b\bigg) =  \sum_{n=0}^\infty R(n)\,q^n.
    \]
We are interested in the situation where $R(0)=1$ and $R(n)\equiv 0 \pmod2$ for $n>0$, i.e., the situation where
$\sum_n R(n) q^n =1$ in the ring of power series $\field{2}[[q]]$. In this case, we say that $A$ and $B$ are
\emph{reciprocals}, and we write $\bar A=B$ and $\bar B=A$. The general problem of this paper is to find the
reciprocals of several special sets $A$, and to draw some conclusions about ``typical'' properties of
reciprocals. We are particularly concerned with the relative density,
    \[
    \delta( \bar A , n) := \frac{| \bar A \cap [0,n]|}{n+1},
    \]
and the density $\delta(\bar A) := \lim_{n\to\infty} \delta(\bar A,n)$ (when the limit exists).

We began studying this problem after reading two articles by Berndt, Yee, and
Zaharescu~\cites{MR2039324,MR1984662}, where bounds on the density of the set
    $
    P_{\text{odd}} := \{n\in \NN \colon p(n) \equiv 1 \pmod{2}\},
    $
with $p(n)$ being the ordinary partition function\footnote{$p(n)$ is the number of ways to write $n$ as a sum of
nonincreasing positive integers. For example, $4=3+1=2+2=2+1+1=1+1+1+1$, so $p(4)=5$.}, are proved. The starting
point for their work is Euler's pentagonal number theorem~\cite{MR1083765}*{Theorem 10.9}, and in particular
that the reciprocal of $P_{\text{odd}}$ is the set $\{ n(3n+1)/2 \colon n\in \ZZ\}$ of pentagonal numbers. Since
the known bounds (see \cite{MR1816213} and \cite{MR1657968} for the currently-best results) on the thickness of
$P_{\text{odd}}$ are so strikingly far from what is believed to be true, we felt that it would be beneficial to
study the ``reciprocal'' notion in a more general setting.

\begin{figure}[t]
\begin{center}
\begin{picture}(360,225)
   \put(0,0){\includegraphics[width=360pt]{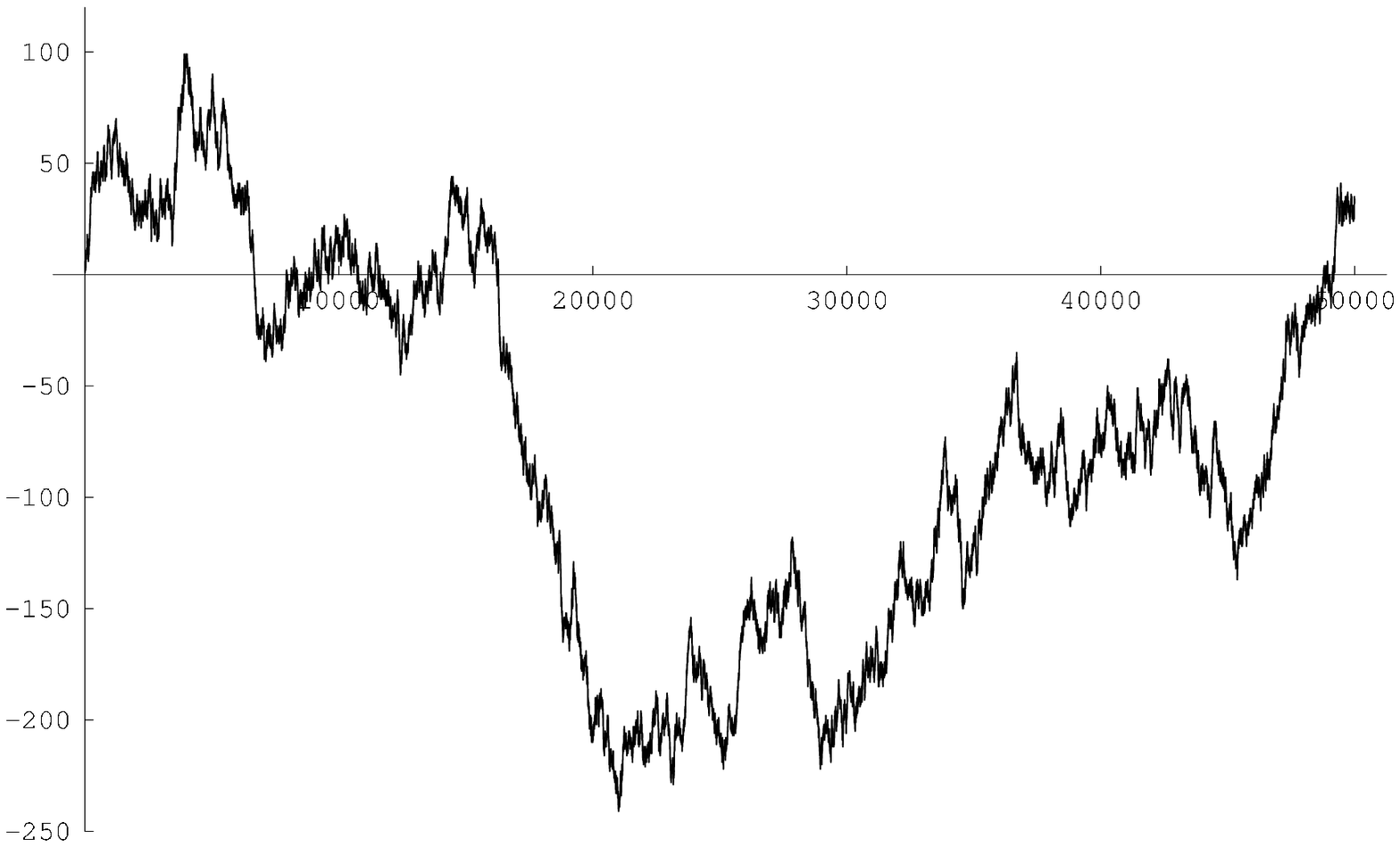}}
\end{picture}
\end{center}
\caption{The points $\big(n,\big| P_{\text{odd}} \cap [0,n]\big|-\big|P_{\text{even}} \cap [0,n]\big| -1\big)=
\big( n, 2\big| P_{\text{odd}} \cap [0,n]\big|-n\big)$\label{fig:walk}}
\end{figure}

$P_{\text{odd}}$ is pictured in Figure~\ref{fig:walk}, where not only does it appear to have density $1/2$, but
the walk defined by $w(n):=2\,\big| P_{\text{odd}}\cap[0,n]\big|-n$ visually appears to be a simple random walk.
See~\cite{MR0227126} for a report of more elaborate statistical tests on the set $P_{\text{odd}}$. We note that
while $p(n)$ appears to be uniformly distributed modulo 2 and 3, it has been known since the time of Ramanujan
to {\em not} be uniformly distributed modulo 5, 7, or 11.

\begin{figure}[t]
\begin{center}
\begin{picture}(380,125)
    \put(0,0){\includegraphics[width=180pt]{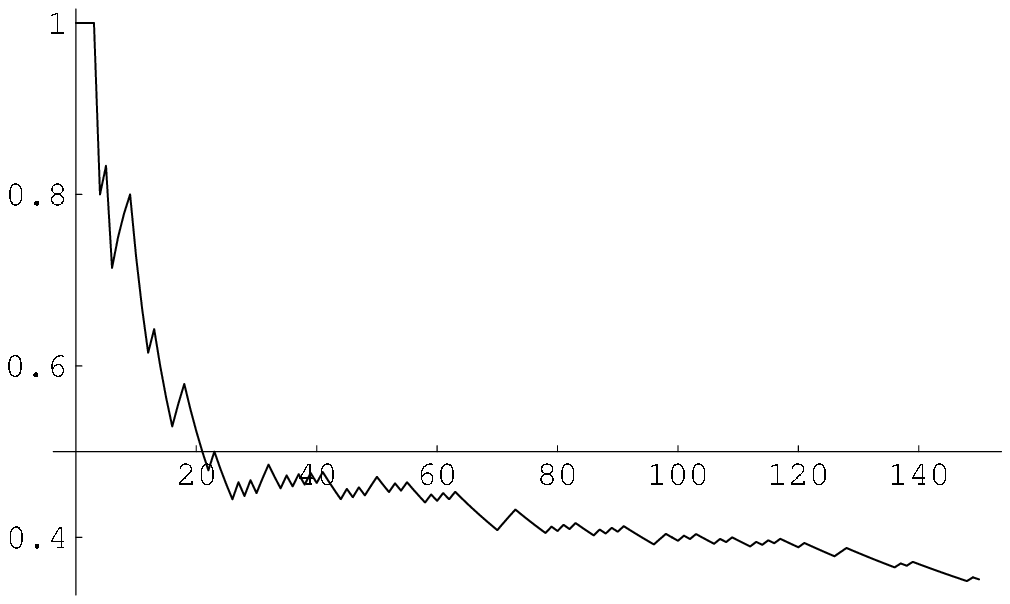}}
    \put(200,0){\includegraphics[width=180pt]{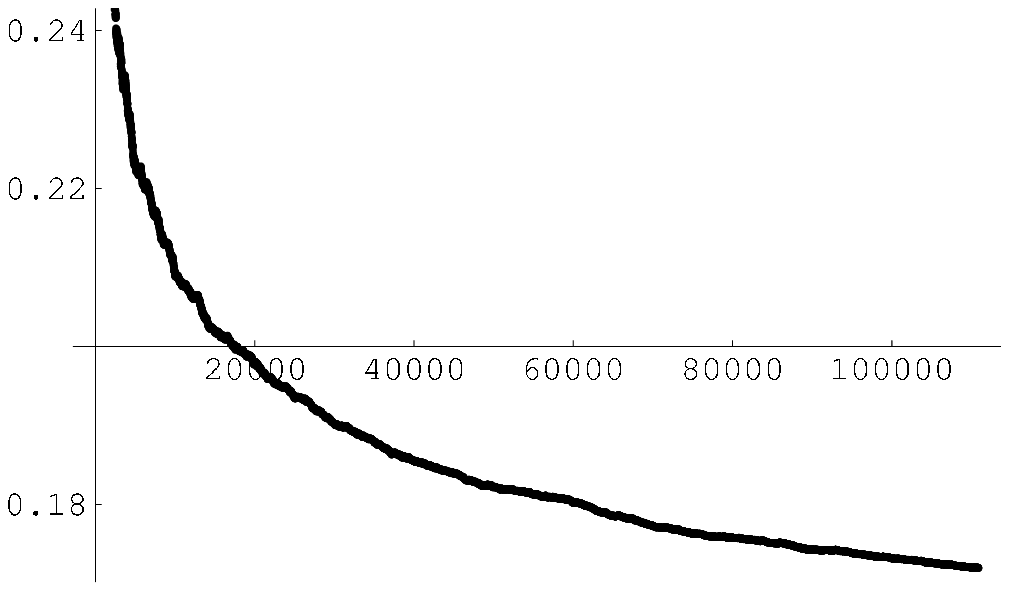}}
\end{picture}
\caption{The relative density of the reciprocal of the set of squares\label{SquarePics}}
\end{center}
\end{figure}

In contrast to that of the pentagonal numbers, the density of the reciprocal of the squares appears to drop off
steadily to 0. Set $S:=\{n^2\colon n\in \NN\}$, with reciprocal $\bar S$.  The relative density of $\bar S$ is
pictured in Figure~\ref{SquarePics}. In Section~\ref{sec:squares}, we prove that the even numbers in $\bar S$
are precisely $\{2n^2\colon n\in\NN\}$, and we characterize the $n\in\bar S$ with $n\equiv 1 \pmod4$ as those
$n$ whose prime factorization has a particular shape. Those $n\in\bar S$ with $n\equiv 3 \pmod 4$ are
characterized in terms of the number of representations of $n$ by certain quadratic forms.

Generalizing the squares and pentagonal numbers, we treat
    \[
    \Theta(c_1,c_2):=\bigg\{c_1 n + c_2 \frac{n(n-1)}2  \colon n\in \ZZ\bigg\}
    \]
for general $c_1$ and $c_2$ in Section~\ref{sec:Theta}. A few interesting special cases are the binomial
coefficients $\Theta(0,1)=\{\tbinom n2 \colon n\in \NN\}$, the squares $\Theta(1,2)$, and the pentagonal numbers
$\Theta(1,3)$.

\begin{cnj}\label{c1c2Conjecture}
The reciprocal of the set $\Theta(c_1,c_2)$, where $0 \leq 2c_1 \le c_2$ and $\gcd(c_1,c_2)=1$, has density $0$
if $c_2 \equiv 2 \pmod{4}$, and otherwise has density $1/2$. More precisely, if $c_2\equiv 2 \pmod 4$, then
    \[\lim_{n\to\infty} \frac{\big|\overline{\Theta(c_1,c_2)}\cap[0,n]\big|}{n/\log n} = C,\]
for some positive constant $C$ depending only on $c_2$. If $c_2\not\equiv 2\pmod4$, then
    \[\limsup_{n\to\infty} \left|\frac{\big|\overline{\Theta(c_1,c_2)}\cap[0,n]\big| -n/2}{\sqrt{n\log\log(n)/2}} \right|
    =1.\]
\end{cnj}

Numerically, it seems that the constant $C$ is 2 if $c_2=2$ or 6, and $C=4$ if $c_2=10$. We lack sufficient data
to guess the other values. The authors believe that the non-effective $c_2\equiv 2 \pmod4$ case might be
provable by showing that the generating function of $\Theta(c_1,c_2)$ is congruent modulo 2 to an integer-weight
modular form, which has almost all of its Fourier coefficients even. This is outside the scope of this paper,
and we leave it as an area for further study.

The $c_2 \not\equiv 2 \pmod4$ case is motivated by the celebrated law of the iterated logarithm. Let $X_1,
X_2,\dots$ be independent random variables taking the values 0 and 1 with probability $1/2$. The law of the
iterated logarithm states that
    \[
    \limsup_{n\to \infty} \left| \frac{\sum_{i=1}^n X_i - n/2}{\sqrt{n\log\log(n)/2}}\right| =1
    \]
with probability 1. What we actually would like to conjecture is that reciprocal of $\Theta(c_1,c_2)$, with
$c_2\not\equiv 2\pmod4$, is statistically indistinguishable from a truly random set with density $1/2$. The
phrase ``statistically indistinguishable'' is too vague, however, so in Conjecture~\ref{c1c2Conjecture} we have
settled for this one specific statistic.

The natural expectation is that, barring some cosmic coincidence or obvious structure, the reciprocal of a set
should have density $1/2$. This is affirmed by the case of a random set, which we handle in detail in
Section~\ref{sec:Random}: let $X_1, X_2,\ldots$ be independent random variables taking the values 0 and 1, with
probabilities bounded away from 0 and 1, and set $F:=\{0\}\cup\{n \colon X_n=1\}$. Theorem~\ref{thm:Random}
states that the reciprocal of $F$ has density $1/2$ with probability 1. This makes the sets whose reciprocals do
not have density $1/2$ the interesting ones.

Our purpose is to identify relevant properties of those sets whose reciprocals have density different from
$1/2$. Specifically, in addition to random sets and $\Theta(c_1,c_2)$, we consider finite sets, the set of
powers of two, and the set of Prouhet-Thue-Morse numbers.
    \begin{itemize}
    \item Finite sets: the reciprocal has a rational density, and appears to
            typically have density slightly below $1/2$. We identify through algebraic
            properties two infinite classes of polynomials, one whose reciprocals have density strictly larger
            than $1/2$, and one whose reciprocals have density at most $1/2$.
    \item Powers of 2: the reciprocal of the thin set $\{0\}\cup\{2^n\colon n\in\NN\}$ is the thin set
            \mbox{$\{2^n-1 \colon n\in \NN\}$.} In particular, we describe the reciprocal of
            $\{0\}\cup\{2^{mn} \colon n\in\NN\}$ for every $m\in\NN$.
    \item Prouhet-Thue-Morse numbers\footnote{$T=\{ 0, 3, 5, 6, 9, 10, 12,15, 17, 18, 20,\dots\}$}: the reciprocal of
            \[
            T := \{n\in\NN\colon \text{the binary expansion of $n$ contains an even number of ``1''s}\}.
            \]
            has density $1/3$. Specifically, we prove that $k\in \overline{T}$ if and only if $k=0$ or
            $(k\pm1)/4$ is an integer whose binary expansion ends in an even number of zeros.
    \end{itemize}
The strongest conjecture that is consistent with our theorems, our experiments, and
Conjecture~\ref{c1c2Conjecture}, is Conjecture~\ref{bigConjecture}.
\begin{cnj}\label{bigConjecture}
If a set contains 0, is not periodic, and is uniformly distributed modulo every power of 2, then its reciprocal
has positive density.
\end{cnj}

We now include a section-by-section agenda for the remainder of the paper.
\begin{description}
    \item[Section \ref{sec:Introduction}:] Motivate and contextualize reciprocals of sets.
    \item[Section \ref{sec:Notation}:] Introduce notation and derive general expressions for reciprocals.
    \item[Section \ref{sec:Random}:] Consider reciprocals of random sets with positive density.
    \item[Section \ref{sec:Polynom}:] Consider reciprocals of finite sets.
    \item[Section \ref{sec:PowersOf2}:] Consider the reciprocal of the powers of 2, and similar sets.
    \item[Section \ref{sec:Theta}:] Consider the reciprocal of $\Theta(c_1,c_2)$, particularly the squares.
    \item[Section \ref{sec:Thue}:] Consider the Prouhet-Thue-Morse sequence.
\end{description}

\section{Notation and General Formulas}\label{sec:Notation}

Throughout this paper, we let
    \begin{equation} \label{eq:fdef}
    \cF(q) = f_0+f_1 q+f_2 q^2 +\cdots \qquad\text{and}\qquad
    \bar{\cF}(q) = \bar f_0+\bar{f_1} q +\bar{f_2} q^2+\cdots
    \end{equation}
be elements of $\field{2}[[q]]$ that satisfy the equation
    \begin{equation} \label{eq:Reciprocals}
    \cF(q)\bar{\cF}(q)=1.
    \end{equation}
In particular, $f_0=\bar f_0 =1$. We define the integer sets $F:=\{n\ge0\colon f_n=1\}$ and $\bar{F}:=\{n\ge 0
\colon \bar f_n=1\}$.

Note that \eqref{eq:Reciprocals} implies (for all $k\ge1$) that $\cF(q^k)\bar{\cF}(q^k)=1$ also. This
corresponds to noting that multiplying everything in $F$ by $k$ has the effect of multiplying everything in
$\bar F$ by $k$. With this in mind, we sometimes make the convenient assumption that $\gcd F =1$.

Our next lemma is a fundamental identity in $\field{2}[[q]]$, and has a number of remarkable consequences. We
use it frequently throughout this paper.
\begin{lem}\label{lem:golden}
The reciprocal of $\cF(q)$ is $\cF(q) \cF(q^2)\cF(q^4)\cF(q^8) \cdots$. That is,
\begin{equation}\label{eq:golden}
1= \cF(q)\,\cdot \, \prod_{k=0}^\infty \cF(q^{2^k}).
\end{equation}
\end{lem}

\begin{proof}
First, notice that both sides of this equation have constant term equal to $1$. Also notice that for any fixed
$n>0$, only finitely many terms of the infinite product affect the coefficient of $q^n$. Thus, the coefficient
of $q^n$ on the right hand side of \eqref{eq:golden} is also the coefficient of $q^n$ in
\begin{eqnarray*}
\cF(q)\,\cdot \, \prod_{k=0}^{\lfloor \log_2 n\rfloor} \cF(q^{2^k}).
\end{eqnarray*}
By the so-called children's binomial theorem\footnote{$(a+b)^2=a^2+b^2 \pmod{2}$}, $\cF(q)\cF(q)= \cF(q^2)$.
Multiplying by $\cF(q^2)$, we see that $\cF(q)\cF(q)\cF(q^2)=\cF(q^2)\cF(q^2)=\cF(q^4)$, and continuing we get
\begin{eqnarray}\label{eq:goldpoly}
    \cF(q) \cF(q) \cF(q^2) \cF(q^4) \cdots \cF(q^{2^{\lfloor \log_2 n\rfloor}})=
    \cF\big(q^{2^{{\lfloor \log_2 n\rfloor}+1}}\big).
\end{eqnarray}
Now notice that since $0<n < 2^{{\lfloor \log_2 n\rfloor}+1}$, the coefficient of $q^n$ on the right hand side
of~\eqref{eq:goldpoly} is $0$, and our result follows.
\end{proof}

We now give a list of recurrences for $\bar f_n$, discuss the usefulness of each, and prove them.

\begin{lem}\label{lem:barfn} If $\cF(q)\bar{\cF}(q)=1$, then $\bar f_0=1$ and
for $n>0$,
    \begin{enumerate}\renewcommand{\theenumi}{\roman{enumi}}
    \item $\displaystyle \bar{f}_n = \sum_{j=1}^{n} f_j \bar{f}_{n-j}$;
    \item $\bar f_n =1$ if and only if $\displaystyle \#\bigg\{(x_0,x_1,\dots)
                                        \colon x_i \in F, n=\sum_{i\ge0} x_i 2^i\bigg\}$ is odd;
    \item $\displaystyle \bar f_n = \sum_{\vec{x}} f_{x_1}f_{x_2}\cdots f_{x_\ell}$,
        where the summation extends over all tuples $\vec{x}=(x_1,\cdots,x_\ell)$ with
        $n=\sum_{i=1}^\ell x_i$ and each $x_i > 0$ ($\ell$ is allowed to vary);
    \item $\displaystyle \bar f_n = \sum_{0\le i<n/4} f_{n-2i}\bar f_i + G(f_1,f_2,\dots,f_{\lfloor n/2\rfloor})$,
    for some function $G$.
    \end{enumerate}
\end{lem}

Lemma~\citebarfnitem{i} is valuable because of its simplicity. For instance, it is immediately apparent from
this recurrence relation that $\bar \cF$ is uniquely defined and always exists (provided $f_0=1$).

In several of the examples we consider, the set $F$ has some special properties modulo a power of 2.
Lemma~\citebarfnitem{ii} facilitates our exploitation of these special properties.

Lemma~\citebarfnitem{iii} is useful because of its symmetry, and because its right-hand side does not expressly
reference the $\bar f$ sequence. As a specific example, let $r(n)$ be the number of ways to write $n$ as a sum
of positive pentagonal numbers (counting order). Then, by Lemma~\citebarfnitem{iii}, $p(n)\equiv r(n)\pmod{2}$.
We also use Lemma~\citebarfnitem{iii}, for example, to prove Lemma~\citebarfnitem{iv}.

If one lets the $f_i$ be independent random variables, then the expression in Lemma~\citebarfnitem{iv} contains
a summation of weakly dependent random variables, and a deterministic function of $f_1,\dots,f_{n/2}$. This
allows us to say something explicit about the distribution of the resulting random variable $\bar f_n$ (see
Theorem~\ref{thm:Random}).

Another remarkable aspect of Lemma~\citebarfnitem{iv} is that $\bar f_n$ does not depend in any way on $f_{n-1},
f_{n-3}, \dots, f_{n-c}$, where $c$ is the largest odd number strictly less than $n/2$. For example,
    \[
    \bar f_{11} = f_{11}+ f_9f_1+f_7f_2+f_7f_1+f_5f_3 + f_4f_3+f_3f_2f_1+f_3f_2+f_2f_1+f_1
    \]
does not depend on $f_{10}, f_8,$ or $f_6$.

\begin{proof}
Comparing the coefficients of $q^n$ on the left- and right-hand sides of equation~\eqref{eq:Reciprocals} yields
$\sum_{j=0}^n f_j \bar{f}_{n-j} = 0$. Lemma~\citebarfnitem{i} is this expression rearranged, using the fact that
$f_0=1$.

Similarly, Lemma~\citebarfnitem{ii} equates the coefficients of $q^n$ on the left- and right-hand sides of
equation~\eqref{eq:golden}, with the right-hand side interpreted as a product in $\ZZ$.

One can prove Lemma~\citebarfnitem{iii} by induction, using Lemma~\citebarfnitem{i} to complete the induction
step. Alternatively, one may simply compare the coefficients of $q^n$ on the left- and right-hand sides of
    \[\bar{\cF}=\frac{1}{\cF}=\frac{1}{1-(\cF-1)}=1+(\cF-1)+(\cF-1)^2+(\cF-1)^3+\cdots,\]
which is valid because we are working over $\field{2}$.

Recall Kummer's result that the multinomial coefficient
    $\displaystyle \binom{m_1+\cdots + m_k}{m_1,m_2,\dots,m_k} = \frac{(m_1+\cdots+ m_k)!}{m_1!m_2!\cdots m_k!}$
is relatively prime to a prime $p$ if and only if $m_1, \ldots, m_k$ can be added in base $p$ without
carrying~\cite{Kummer}. We are working with $p=2$, so our condition is: $\tbinom{m_1+\cdots+
m_k}{m_1,\dots,m_k}$ is odd if and only if no two of the binary expansions of $m_1,\dots,m_k$ have a ``1'' in
the same position. We call such a list of positive integers $m_1,\dots,m_k$ {\em non-overlapping}.

Let $\pi(n)$ be the set of partitions of $n$ whose distinct parts $x_1,\dots,x_k$ have non-overlapping
multiplicities $m_1,\dots,m_k$. Continuing from \citebarfnitem{iii}, we have
    \begin{align*}
    \bar f_n &= \sum_{\substack{x_1+\cdots+x_\ell=n \\x_i>0}} f_{x_1}f_{x_2}\cdots f_{x_\ell} \\
             &= \sum_{\substack{m_1a_1+\cdots +m_ka_k=n\\ a_1>\cdots>a_k>0\\m_i>0}}\binom{m_1+\cdots+m_k}{m_1,\dots,m_k}
                        f_{a_1}^{m_1} f_{a_2}^{m_2} \cdots f_{a_k}^{m_k}\\
             &= \sum_{\pi(n)} f_{a_1}^{m_1} f_{a_2}^{m_2} \cdots f_{a_k}^{m_k}\\
             &= \sum_{\substack{\pi(n)\\a_1>n/2}} f_{a_1}^{m_1} f_{a_2}^{m_2} \cdots f_{a_k}^{m_k} +
                    \sum_{\substack{\pi(n)\\a_1\le n/2}} f_{a_1}^{m_1} f_{a_2}^{m_2} \cdots f_{a_k}^{m_k}
    \end{align*}
If $a_1>n/2$, then it must have multiplicity $m_1=1$, and if $m_1,\dots,m_k$ are non-overlapping then the other
$m_i$ are even:
    \[ f_{a_1}^{m_1} f_{a_2}^{m_2} \cdots f_{a_k}^{m_k} = f_{a_1} \left(f_{a_2}^{m_2/2}\cdots
    f_{a_k}^{m_k/2}\right)^2.\]
This implies that $n-a_1$ is even, and $a_2\frac{m_2}2 + \dots + a_k \frac{m_k}2$ is a partition of $(n-a_1)/2$.
Setting $2i=n-a_1$, we get
    \[
    \bar f_n = \sum_{\pi(n)} f_{a_1}^{m_1} f_{a_2}^{m_2} \cdots f_{a_k}^{m_k}
        = \sum_{0\le i < n/4} f_{n-2i} \sum_{\pi(i)} f_{a_1}^{m_1}\cdots f_{a_k}^{m_k}
            + \sum_{\substack{\pi(n)\\a_i\le n/2}}f_{a_1}^{m_1} \cdots f_{a_k}^{m_k}.\]
Using Lemma~\citebarfnitem{iii} again, this becomes
    \[
    \bar f_n =\sum_{0\le i < n/4} f_{n-2i}\bar f_i + G(f_1,\ldots,f_{\lfloor n/2 \rfloor})
    \]
for a specific function $G$.
\end{proof}

\section{Random power series}\label{sec:Random}

     In this section, we consider the reciprocal of a
random power series in $\field{2}[[q]]$. The results of this section are strong evidence that the density of
$\cF$ plays little to no role in determining the density of $\bar\cF$, and that unless the coefficients of $\cF$
have some structure, the density of $\bar\cF$ is $1/2$.

Recall that a Bernoulli variable is a random variable that takes values in $\{0,1\}$.

\begin{thm}\label{thm:Random}
Suppose that $f_1, f_2, \dots$ are independent Bernoulli variables, with
    \[\inf_n \min\{\Prob{f_n=0},\Prob{f_n=1}\}>0.\]
Then $\delta(\bar{\cF})=1/2$ with probability 1.
\end{thm}

We need the following two lemmas.

\begin{lem}[L\'{e}vy's Borel-Cantelli lemma] \label{lem:LBCL}
Let $E_1, E_2, \dots, $ be events, and define $Z_n:=\sum_{k=1}^n
I_{E_k}$, the random variable that records the number of $E_1, E_2, \dots, E_n$ that occur. Define
    \[\xi_k:=\Prob{E_k \mid E_1, E_2, \dots, E_{k-1}}.\]
If $\sum_{k=1}^\infty \xi_k$ diverges, then $Z_n$ is asymptotically equal to $\sum_{k=1}^n \xi_k$ with
probability 1.
\end{lem}

For a proof of L\'{e}vy's Borel-Cantelli lemma, we refer the reader to \cite{MR1155402}*{Sec 12.15}.

\begin{lem}[Binary Central Limit Theorem] \label{lem:BCLT}
Let $X_i$ be 0 with probability $\gamma_i$ and 1 with probability $1-\gamma_i$, and suppose that $X_1,
X_2,\dots$ are independent. Then, as $n\to\infty$,
    \[\Prob{\sum_{i=1}^n X_i \equiv 0 \pmod 2} \to \frac12\]
if and only if some $\gamma_i=1/2$ or $\sum_{i=1}^n \min\{\gamma_i,1-\gamma_i\}$ diverges.
\end{lem}

\begin{proof} Let $S_n:=\sum_{i=1}^n X_i$, and define $p_i$ by $\Prob{S_n\equiv 0 \pmod{2}}=p_i$.
Clearly $S_{n}$ is even if and only if $S_{n-1}$ and $X_{n}$ are both even or both odd:
    \[
    p_{n} = p_{n-1} \gamma_{n} +(1-p_{n-1}) (1-\gamma_{n}).
    \]
Clearly $2p_1-1=2\gamma_1-1$, and
    \[
    2p_{n}-1    =   2\left(p_{n-1} \gamma_{n} +(1-p_{n-1}) (1-\gamma_{n})\right)-1
                =  (2p_{n-1}-1)(2\gamma_{n}-1),
    \]
which provides the base case and inductive step for the equality
    \[2p_n-1 =  \prod_{i=1}^{n} (2\gamma_i-1).\]
By the standard results for infinite products, we now see that $2p_n-1\to0$ if and only if $2\gamma_i-1=0$ for
some $i$ or $\sum_{i=1}^n \min\{\gamma_i,1-\gamma_i\}$ diverges.
\end{proof}

\begin{proof}[{Proof of Theorem~\ref{thm:Random}}]
We begin with some notation:
\begin{align*}
    \alpha_n &:= \Prob{f_n=0}, \\
    \beta_n &:=\min\{\alpha_n,1-\alpha_n\},\\
    \beta &:=\inf_{n\to\infty} \beta_n,
\end{align*}
and note that $0<\beta \le 1/2$. We will show first that $\Prob{\bar{f}_n=0} \rightarrow 1/2$ as $n \rightarrow
\infty$, and then will show that $\delta(\bar F)=1/2$ with probability 1.

Lemma~\citebarfnitem{i} says that ${\bar f_n} = {f_n+ \sum_{j=1}^{n-1}f_j\bar f_{n-j}}$, whence
    \[
        \Prob{\bar f_n=0}=\Prob{f_n=0}\Prob{\sum_{j=1}^{n-1}f_j\bar f_{n-j}=0}
                            + \Prob{f_n=1}\Prob{\sum_{j=1}^{n-1}f_j\bar f_{n-j}=1}
    \]
is a weighted average of $\Prob{f_n=0}=\alpha_n$ and $\Prob{f_n=1}=1-\alpha_n$. Consequently, $\Prob{\bar
f_n=0}\geq \beta_n \geq \beta$ and $\Prob{\bar f_n=1} \ge \beta_n\ge \beta$.

Set $B_n:=\{i \colon 0\le i < n/4, \bar f_i=1\}$, and set $G_n:=G(f_1,\dots,f_{\lfloor n/2\rfloor})$, where $G$
is the function from Lemma~\citebarfnitem{iv}. We have, from Lemma~\citebarfnitem{iv},
    \[\bar f_n = \sum_{i\in B_n} f_{n-2i} + G_n.\]
From the previous paragraph, we know that $\Expect{|B_n|}$ is at least $\sum_{i=0}^{\lfloor n/4 \rfloor} \beta_i
\geq \beta\lfloor n/4 \rfloor$. In particular, a routine calculation shows that $|B_n|\to\infty$ with
probability 1. Thus, $\Prob{|B_n|>K_n}\to 1$ if $K_n$ goes to infinity sufficiently slowly. We have
    \begin{multline*}
    \Prob{\bar f_n=0} =\\
                \Prob{\bar f_n=0 \,\big|\, |B_n|>K_n}\Prob{|B_n|>K_n}
                +\Prob{\bar f_n=0 \,\big|\, |B_n|\le K_n}\Prob{|B_n|\le K_n},
    \end{multline*}
which for large $n$ becomes $\Prob{\bar f_n=0}=\Prob{\bar f_n=0 \,\big|\, |B_n|>K_n}$.

We now observe that $\bar f_n=0$ if and only if $G_n=\sum_{i\in B_n} f_{n-2i}$ (call this sum $\sigma_n$), so
that
    \begin{multline*}
    \Prob{\bar f_n=0 \,\big|\, |B_n|>K_n} = \\
        \hspace{1cm}\Prob{G_n=0 \,\big|\, |B_n|>K_n, \sigma_n=0}\Prob{\sigma_n=0 \,\big|\, |B_n|>K_n} +\\
        \Prob{G_n=1\,\big|\, |B_n|>K_n, \sigma_n=0}\Prob{\sigma_n=1 \,\big|\, |B_n|>K_n}.
    \end{multline*}
This is a weighted average of $\Prob{\sigma_n=0 \,\big|\, |B_n|>K_n}$ and $\Prob{\sigma_n=1 \,\big|\,
|B_n|>K_n}$, both of which go to $1/2$ as $n\to\infty$ by the Binary Central Limit Theorem. Thus,
    \[
    \Prob{\bar f_n=0} \approx \Prob{\bar f_n=0 \,\big|\, |B_n|>K_n} \approx \frac 12
    \]
with each of the ``$\approx$'' becoming ``$=$'' as $n\to \infty$.

Now that we have shown that $\Prob{\bar f_n=0}\to1/2$, we know that $\Expect{\delta(\bar F,n)} \to 1/2$, but
this does not imply that $\delta(\bar F,n)\to1/2$ ever, much less with probability 1. This last step again
requires the at-least-weak independence of $\bar f_n$ from $\bar f_1, \dots, \bar f_{n-1}$, and the
technicalities are handled for us by L\'{e}vy's Borel-Cantelli lemma.

Let $E_k$ be the event $\{\bar f_k =0\}$, and set $\xi_k:=\Prob{\bar f_k=0 \mid \bar f_1, \bar f_2, \dots, \bar
f_n}$. By the comment above, $0<\beta\le \xi_k$, so $\sum_{k=1}^\infty \xi_k=\infty$. Thus, by
Lemma~\ref{lem:LBCL},
    \begin{equation}\label{eq:asymptdef}
    \lim_{n\to\infty} \frac{\delta(\bar F,n)}{\frac 1n \sum_{k=1}^n \xi_k} = 1
    \end{equation}
with probability 1. For every $\epsilon>0$ there is an $n_0$ such that for all $n>n_0$
    \[
    (1-\epsilon) \frac 1n \sum_{k=1}^n \xi_k \leq \delta(\bar F,n) \le (1+\epsilon) \frac 1n \sum_{k=1}^n \xi_k.
    \]
These upper and lower bounds on $\delta(\bar F,n)$ are non-random, so we may take expectations (for large $n$)
to get
    \[
    (1-\epsilon) \frac 1n \sum_{k=1}^n \xi_k \leq \frac 12 \le (1+\epsilon) \frac 1n \sum_{k=1}^n \xi_k,
    \]
where we have used the linearity of expectation and the previously proved $\Expect{\bar f_n=0}=\Prob{\bar f_n=0}
\to 1/2$. This implies that $\frac 1n\sum_{k=1}^n \xi_k \to 1/2$ also. Consequently, \eqref{eq:asymptdef} now
implies that
    \[
    \delta(\bar F) := \lim_{n\to\infty} {\delta(\bar F,n)} = \frac 12
    \]
with probability 1.
\end{proof}

%
%

\section{Polynomials }\label{sec:Polynom}

In this section, we study the reciprocals of polynomials in $\field{2}[q]$. The coefficients of such a
reciprocal are periodic (see Proposition~\ref{lem:lsr} below), and so the reciprocal has rational density. We
also give some indication of how the densities of reciprocals of polynomials are distributed, beginning in
Subsection~\ref{sec:orderanddensity}. In Subsection~\ref{sec:dbca}, we use the theory of de Bruijn cycles to
exhibit an infinite family of polynomials whose reciprocals have densities strictly larger than $1/2$; in
Subsection~\ref{sec:companions}, we show that if two polynomials have product $1+q^D$ ($D\ge4$), then at least
one of them has a reciprocal with density at most $1/2$. In Subsection~\ref{sec:eps}, we show that the
reciprocal of an eventually periodic set\footnote{More precisely, a set whose indicator function is eventually
periodic.} containing 0 is an eventually periodic set containing 0.

Let $\sum_{i=0}^\infty b_i 2^i$ be the binary expansion of $n$; we define the polynomial
$\cP_n(q):=\sum_{i=0}^\infty b_i q^i \in \field{2}[q]$. Clearly this indexes all polynomials, and the invertible
polynomials are precisely those with $n$ odd. For a polynomial $\cQ\in \field{2}[q]$, we let $\hat\cQ$ be the
same polynomial with coefficients (all 0 or 1) in $\ZZ$. For instance, $\cQ=\cP_{\hat \cQ(2)}$ for every
polynomial $\cQ$.

We denote by $\ell(\cP)$ the length of the polynomial, i.e., $\ell(\cP)=\hat{\cP}(1)$, and by $\deg(\cP)$ the
degree of the polynomial. Also, $\ord(\cP)$ is the least positive $D$ such that $\cP$ divides $1+q^D$. It is not
immediately obvious that $\ord(\cP)$ is well defined: it is for invertible $\cP$, and this is the content of
Proposition~\ref{lem:ordfinite} below. For each polynomial $\cP$, we define $\cP^\ast$ by
$\cP\cP^\ast=1+q^{\ord{\cP}}$. We shall see that the properties of $\cP$ and $\cP^\ast$ are intimately related
(Propositions~\ref{lem:densityOfF} and~\ref{prop:1/2}).

If $\ord(\cP) = 2^{\deg(\cP)}-1$, then $\cP$ is called {\it primitive}, and $\field{2}[q]/(\cF)$ is isomorphic
to $\field{2^{\deg(\cP)}}$, with multiplicative generator $q$. All primitive polynomials are irreducible, but
not vice versa; for example $1+q^3+q^6$ and $1+q+q^2+q^3+q^4$ are irreducible but {\em not} primitive.

Figure~\ref{polydata} tabulates properties of $\cP_n$ for odd $n<256$, including factorizations, $\cP_n^\ast$,
and densities of reciprocals.

\begin{figure}[p!]
\begin{center}
 {\tiny
\begin{tabular}{ccccc|ccccc}
 $n$ & $\hat{\cP}^\ast(2)$ & $D$ & Factors & $\delta(\bar \cP_n)$ & $n$ &  $D$ & Factors & $\delta(\bar \cP_n)$ \\ \hline
    1 & 1 & 1 & prim & 0 & 129 &  7 & $ 3 \cdot 11 \cdot 13 $ & 1/7 \\
    3 & 1 & 1 & prim & 1 & 131  & 127 &     prim & 64/127 \\
    5 & 1 & 2 & $ 3^2 $ & 1/2 & 133  & 93 & $ 7 \cdot     55 $ & 46/93 \\
    7 & 3 & 3 & prim & 2/3 & 135  & 60 & $ 3^3 \cdot 25 $ & 1/2     \\
    9 & 1 & 3 & $ 3 \cdot 7 $ & 1/3 & 137   & 127 & prim & 64/127 \\
    11 & 23 & 7 & prim & 4/7 & 139  & 15 & $ 3 \cdot 7 \cdot 19 $ & 1/3 \\
    13 & 29 & 7 & prim & 4/7 & 141  & 62 & $ 3^2 \cdot 41 $ &     1/2 \\
    15 & 3 & 4 & $ 3^3 $ & 1/2 & 143  &     127 & prim & 64/127 \\
    17 & 1 & 4 & $ 3^4 $ & 1/4 & 145  &     127 & prim & 64/127 \\
    19 & 2479 & 15 & prim & 8/15 & 147  & 62 & $ 3^2 \cdot 47     $ & 1/2 \\
    21 & 5 & 6 & $ 7^2 $ & 1/3 & 149  & 63 & $ 3 \cdot 115 $ &     31/63 \\
    23 & 11 & 7 & $ 3 \cdot 13 $ & 3/7 & 151  & 42 & $ 7^2 \cdot 11     $ & 10/21 \\
    25 & 3929 & 15 & prim & 8/15 & 153  & 24 & $ 3^5 \cdot 7 $ & 1/2 \\
    27 & 7 & 6 & $ 3^2 \cdot 7 $ & 1/2 & 155  & 35 & prim & 18/35 \\
    29 & 13 & 7 & $ 3 \cdot 11 $ & 3/7 & 157  & 127 & prim & 64/127 \\
    31 & 3 & 5 & irr & 2/5 & 159  & 21 & $ 3 \cdot 117 $ & 3/7 \\
    33 & 1 & 5 & $ 3 \cdot 31 $ & 1/5 & 161  & 93 &     $ 7 \cdot 59 $ & 46/93 \\
    35 & 72031 & 21 & $ 7 \cdot 13 $ & 10/21 & 163  & 63 & $ 3     \cdot 97 $ & 31/63 \\
    37 & 78898037 & 31 & prim & 16/31 & 165  & 20 & $ 3^3 \cdot 31 $ & 1/2     \\
    39 & 635 & 14 & $ 3^2 \cdot 11 $ & 1/2 & 167  & 127 & prim & 64/127 \\
    41 & 91635305 & 31 & prim & 16/31 & 169  & 63 & $ 3 \cdot     103 $ & 31/63 \\
    43 & 1335 & 15 & $ 3 \cdot 25 $ & 7/15 & 171  & 127 & prim & 64/127 \\
    45 & 189 & 12 & $ 3^3 \cdot 7 $ & 1/2 & 173      & 105 & $ 11 \cdot 19 $ & 52/105 \\
    47 & 94957459 & 31 & prim & 16/31 & 175  & 42 & $ 3^2 \cdot 7     \cdot 13 $ & 1/2 \\
    49 & 128305 & 21 & $ 7 \cdot 11 $ & 10/21 & 177  & 62 & $     3^2 \cdot 37 $ & 1/2 \\
    51 & 15 & 8 & $ 3^5 $ & 1/2 & 179  & 93 & $ 7     \cdot 61 $ & 46/93 \\
    53 & 1893 & 15 & $ 3 \cdot 19 $ & 7/15 & 181   & 105 & $ 13 \cdot 25 $ & 52/105 \\
    55 & 121098539 & 31 & prim & 16/31 & 183  & 63 & $ 3     \cdot 109 $ & 31/63 \\
    57 & 889 & 14 & $ 3^2 \cdot 13 $ & 1/2 & 185  & 127 & prim & 64/127 \\
    59 & 111435623 & 31 & prim & 16/31 & 187  & 28 & $ 3^4 \cdot 11 $ &     13/28 \\
    61 & 105887917 & 31 & prim & 16/31 & 189  & 12 & $ 3 \cdot 7^3 $ & 1/3 \\
    63 & 3 & 6 & $ 3 \cdot 7^2 $ & 1/3 & 191  & 127 & prim & 64/127 \\
    65 & 1 & 6 & $ 3^2 \cdot 7^2 $ & 1/6 & 193  & 127 & prim & 64/127 \\
    67 & 151054908502416063 & 63 & prim & 32/63 & 195  & 12 & $ 3^3 \cdot 7^2     $ & 1/2 \\
    69 & 277 & 14 & $ 11^2 $ & 2/7 & 197  & 63 & $ 3 \cdot 67     $ & 31/63 \\
    71 & 37394331 & 31 & $ 3 \cdot 61 $ & 15/31 & 199  & 105 & $ 13 \cdot 19 $ & 52/105 \\
    73 & 9 & 9 & irr & 2/9 & 201  & 62 & $ 3^2 \cdot 61 $ &     1/2 \\
    75 & 4865751 & 28 & $ 3^3 \cdot 13 $ & 1/2 & 203  & 127 & prim & 64/127 \\
    77 & 40094429 & 31 & $ 3 \cdot 59 $ & 15/31 & 205  & 93 & $ 7 \cdot 47 $ & 46/93 \\
    79 & 627 & 15 & $ 7 \cdot 25 $ & 2/5 & 207  & 14 & $ 3 \cdot 11^2 $ &     3/7 \\
    81 & 337 & 14 & $ 13^2 $ & 2/7 & 209  & 15 & $ 3 \cdot 7 \cdot 25 $ &     1/3 \\
    83 & 44271 & 21 & prim & 11/21 & 211   & 127 & prim & 64/127 \\
    85 & 5 & 8 & $ 3^6 $ & 1/4 & 213  &     127 & prim & 64/127 \\
    87 & 42187 & 21 & irr & 8/21 & 215  & 62 & $ 3^2 \cdot 59     $ & 1/2 \\
    89 & 49106713 & 31 & $ 3 \cdot 55 $ & 15/31 & 217  & 35 & prim &     18/35 \\
    91 & 215232491192501383 & 63 & prim & 32/63 & 219  & 9 & $ 3 \cdot 73 $ &     1/3 \\
    93 & 717 & 15 & $ 7 \cdot 31 $ & 2/5 & 221  & 28 & $ 3^4 \cdot 13 $     & 13/28 \\
    95 & 24018211 & 30 & $ 3^2 \cdot 19 $ & 1/2 & 223  & 93 & $ 7 \cdot 41 $ & 46/93 \\
    97 & 285247320157033569 & 63 & prim & 32/63 & 225  & 60 &     $ 3^3 \cdot 19 $ & 1/2 \\
    99 & 31 & 10 & $ 3^2 \cdot 31 $ & 1/2 & 227   & 105 & $ 11 \cdot 25 $ & 52/105 \\
    101 & 63285 & 21 & prim & 11/21 & 229    & 127 & prim & 64/127 \\
    103 & 272840796136989499 & 63 & prim & 32/63 & 231  & 15 & $ 3 \cdot 7     \cdot 31 $ & 7/15 \\
    105 & 7716393 & 28 & $ 3^3 \cdot 11 $ & 1/2 & 233  & 42 & $ 7^2     \cdot 13 $ & 10/21 \\
    107 & 119 & 12 & $ 7^3 $ & 1/2 & 235 & 62 & $ 3^2 \cdot     55 $ & 1/2 \\
    109 & 253483157574931709 & 63 & prim & 32/63 & 237  & 63     & $ 3 \cdot 91 $ & 31/63 \\
    111 & 59858643 & 31 & $ 3 \cdot 37 $ & 15/31 & 239  & 127 & prim & 64/127 \\
    113 & 57124209 & 31 & $ 3 \cdot 47 $ & 15/31 & 241  & 127 & prim & 64/127 \\
    115 & 248574834945763919 & 63 & prim & 32/63 & 243  & 14 & $ 3 \cdot     13^2 $ & 3/7 \\
    117 & 54053 & 21 & irr & 8/21 & 245  & 42 & $ 3^2 \cdot 7 \cdot     11 $ & 1/2 \\
    119 & 107 & 12 & $ 3^4 \cdot 7 $ & 5/12 & 247  & 127 & prim & 64/127 \\
    121 & 825 & 15 & $ 7 \cdot 19 $ & 2/5 & 249  & 21 & $ 3 \cdot 87 $ &     3/7 \\
    123 & 53340711 & 31 & $ 3 \cdot 41 $ & 15/31 & 251  & 93 & $ 7 \cdot 37 $ & 46/93 \\
    125 & 25787629 & 30 & $ 3^2 \cdot 25 $ & 1/2 & 253  & 127 & prim & 64/127 \\
    127 & 3 & 7 & $ 11 \cdot 13 $ & 2/7 & 255 & 8 & $ 3^7 $ & 1/4 \\
\end{tabular}
 }
 \caption{
 Properties of $\cP_n$ for odd $n<256$.
 The ``Factors'' column records whether $\cP_n$ is
 reducible, irreducible or primitive.
 If $\cP_n$ is reducible,
 then the factors column evaluates the factors of $\cP_n$ at $q=2$. For example, the factors of
 $P_{245}$ are given as $ 3^2 \cdot 7 \cdot 11 $, whence $P_{245}=P_{3}^2 P_7 P_{11}$.\label{polydata}}
\end{center}
\end{figure}

\begin{figure}[ht]
\begin{center}
\begin{picture}(400,145)
    \put(26,12){\includegraphics[width=5in]{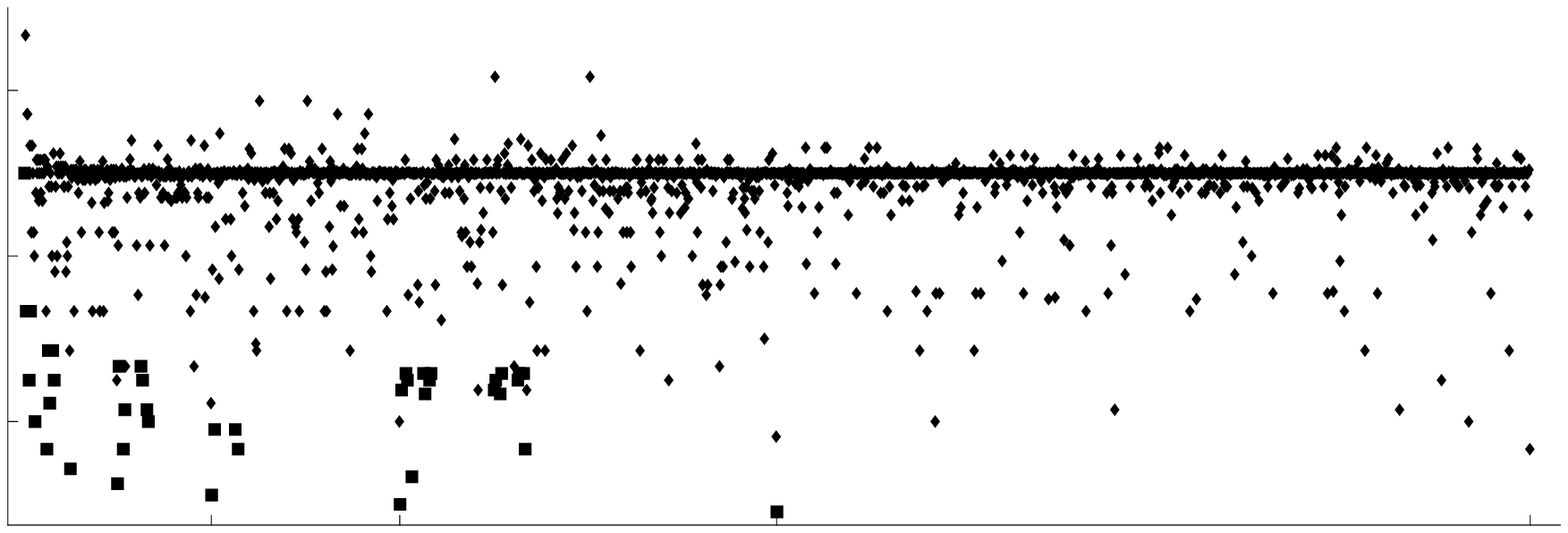}}
    \put(65,0){512}
    \put(104,0){1024}
    \put(193,0){2048}
    \put(364,0){4096}
    \put(10,35){$0.2$}
    \put(10,73){$0.4$}
    \put(10,111){$0.6$}
    \put(390,10){$n$}
    \put(24,136){$\delta$}
\end{picture}
\end{center}
\caption{The points $\big(n,\delta(\bar \cP_n)\big)$ with $n$ odd, except $(1,0)$ and
$(3,1)$\label{PolyInversePic}}
\end{figure}

In Figure~\ref{PolyInversePic}, we plot the points $(n,\delta(\bar\cP_n))$ for odd $n$ less than $2^{12}$. We
note that $\delta(\bar\cP_n)$ tends to be near $1/2$, but is biased toward being below $1/2$. This is also
suggested, but not proven, by Proposition~\ref{prop:1/2} below. In Proposition~\ref{primitivedensity}, we give
an algebraically-described infinite set of $n$ such that $\delta(\bar \cP_n)> 1/2$. Note that $\delta(\bar
\cP_n(q^k))=\frac 1k \delta(\bar\cP_n(q))$, i.e., if $\cP_n$ is a polynomial in $q^2$, $q^3$, etc, then its
density is {\em a priori} less than $1/2$, $1/3$, etc. These points have been plotted with squares.

\begin{figure}[ht]
\begin{center}
\begin{picture}(420,220)
    \put(0,0){\includegraphics[width=5.5in]{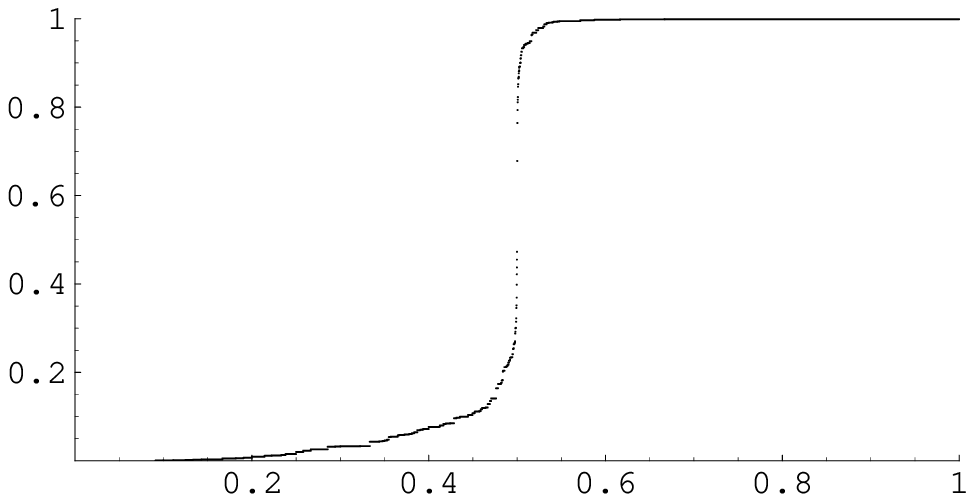}}
    \put(400,10){$x$}
    \put(14,204){$\Delta(x)$}
\end{picture}
\end{center}
\caption{The distribution function $\Delta(x):=2^{-11}\,\cdot\,\#\{n\colon 1\le n \le 2^{12}, n\text{ odd},
\delta(\bar\cP_n)\le x \}$\label{fig:DensityDistribution}}
\end{figure}

In Figure~\ref{fig:DensityDistribution}, we plot the empirical distribution function of $\delta(\bar \cP_n)$.
The large discontinuities near $1/2$ mean that these densities occur with large frequency (fully 421 of the 2048
polynomials $\cP_1,\cP_3,\dots,\cP_{4095}$ have reciprocals with density {\em exactly} $1/2$). Again visible in
Figure~\ref{fig:DensityDistribution} is the preference of $\cP$ to have reciprocal with density less than $1/2$.
The most interesting issue raised in this section, which remains unanswered, is to describe the set
    \[
    \big\{ \delta(\bar \cP) \colon \cP\text{ is a polynomial}\big\}.
    \]
For example, is there an $n$ with $\delta(\bar \cP_n)=3/4$?

\subsection{Order and Density}\label{sec:orderanddensity}

Our first proposition demonstrates that $\ord(\cP)$ is well-defined, and our next proposition shows the
connection between $\delta(\bar \cP)$, $\ord(\cP)$, and $\cP^\ast$.

\begin{prop}\label{lem:ordfinite}
If $\cP$ is a polynomial, then $\ord(\cP)$ is finite.\footnote{Actually, the proof can be refined to show that
$\ord(\cP) \mid 2^{\deg(\cP)}-1$ if $\cP$ is irreducible, and otherwise $\ord(\cP)=2^i
\lcm\{\ord(\cV_1),\dots,\ord(\cV_k)\}$ for some $1\le 2^i \le k$, where $\cP=\cV_1\cdots \cV_k$.}
\end{prop}

\begin{proof}
Let $\cV_1, \dots, \cV_k$ be the irreducible factors of $\cP$, and let $d_i$ be the multiplicative order of $q$
in the field $\field{2}[q]/(\cV_i)$. In particular, $1+q^{x d_i}$ is a multiple of $\cV_i$ for each $x\in\NN$.
Set $L:=\lcm\{d_1,\ldots, d_k\}$ and define $\cV_i^\ast$ by $\cV_1\cV_1^\ast=1+q^{L}$ and for $1<i\le k$ by
$\cV_i \cV_i^\ast = 1+q^{2^{i-2}L}$. Now
    \begin{multline*}
    \cP \, \cdot \, \prod_{i=1}^k \cV_i^\ast = (1+q^L)(1+q^L)(1+q^{2L})\cdots (1+q^{2^{k-2}L}) \\
        = (1+q^{2L})(1+q^{2L})(1+q^{4L}) \cdots (1+q^{2^{k-2}L})
        = (1+q^{2^{k-1}L}),
    \end{multline*}
by repeated use of the children's binomial theorem.
\end{proof}

We emphasize that, given $\cP$ and the equality $\cP\cF=1+q^D$ for some $\cF$, $D$ is not uniquely determined.
For example, $\cP(q)\cP(q)\cF(q^2)=1+q^{2D}$. Nor does the proof given above always provide the minimal $D$.

\begin{prop}\label{lem:densityOfF}
$\displaystyle \delta(\bar \cP) = {\ell(\cP^\ast)}/{\ord(\cP)}$.
\end{prop}

\begin{proof}
Since $\cP \, \frac{\cP^\ast}{1+q^{\ord{\cP}}}=1$, we see that the reciprocal of $\cP$ is periodic with period
$\ord{\cP}$ (although this may not be the minimal period), and in each period has density
$\ell(\cP^\ast)/\ord{\cP}$.
\end{proof}

\subsection{de Bruijn cycle algebra}\label{sec:dbca}

Our next proposition shows that the reciprocal of a polynomial is a special case of a linear-shift register.
Fortunately, there is an enormous literature on linear-shift registers (see \cite{ShiftRegisterBook}, for
example).
\begin{prop}\label{lem:lsr}
If $\cF$ is a polynomial with degree $d$, then (letting $\bar f_{j}=0$ for negative $j$)
    \begin{equation}\label{equ:recurrence}
    \bar f_n= \sum_{j=1}^d f_j \bar f_{n-j}.
    \end{equation}
Alternatively, $\bar f_n$ is the constant term of $q^{-n} \bmod{\cF}$.
\end{prop}

\begin{proof}
Since $f_j = 0$ for all $j>d$, the recurrence~\eqref{equ:recurrence} follows immediately from
Lemma~\citebarfnitem{i}.

Let $M$ be the matrix whose $k^\textrm{th}$ row is the elementary vector supported in coordinate $k+1$, for $k =
1, \ldots, d-1$, and whose last row is the vector $(f_0,\ldots,f_{d-1})$, i.e., $M$ is the companion matrix of
$\cF$. Write $c_n$ for the constant coefficient of $q^{n} \bmod \cF$. We claim that
    \begin{equation} \label{linearrecurrence}
    M \begin{pmatrix} c_{s} \\ \vdots \\ c_{s+d-1} \end{pmatrix} =
    \begin{pmatrix} c_{s+1} \\ \vdots \\ c_{s+d} \end{pmatrix}.
    \end{equation}
To see this, let $Y_k$ denote scalar projection of elements of $\field{2}[q]/(\cF)$ onto $q^k$, and let $X$
denote multiplication by $q$ in $\field{2}[q]/(\cF)$. Both of these maps are linear, and it is easy to see that
    $
    Y_k = Y_{k-1} X^{-1} + f_k Y_0,
    $
for $1 \leq k \leq d-1$. Therefore,
    \begin{align*}
    Y_0 & = Y_{d-1} X^{-1} \\
    & = Y_{d-2} X^{-2} + f_{d-1} Y_0 X^{-1} \\
    & = Y_{d-3} X^{-3} + f_{d-2} Y_0 X^{-2} + f_{d-1} Y_0 X^{-1} \\
    & \,\,\,\, \vdots \\
    & = \sum_{j=0}^{d-1} f_{j} Y_0 X^{d-j}.
    \end{align*}
Applying $Y_0$ to $q^{s+d}$ yields $c_{s+d}=\sum_{j=0}^{d-1} f_{j} c_{s+j}$, which implies
\eqref{linearrecurrence}. Set $a_n:=c_{-n}$ (define $c$ on negative subscripts by using the recurrence). Thus,
the sequences $(a_{n})$ and $(\bar f_n)$ satisfy the same recurrence, with initial conditions $a_0=\bar f_0=1,
c_{-i}=a_{i}=\bar f_{i}=0$ (for $-d<i<0$).
\end{proof}

Our next proposition computes the density of the reciprocal of every primitive polynomial, and thereby produces
an infinite family of polynomials whose reciprocals have density greater than $1/2$.

\begin{prop} \label{primitivedensity}
If $\cP$ is a primitive polynomial with degree $d$, then $\displaystyle \delta(\bar\cP) =
\frac{2^{d-1}}{2^d-1}$.
\end{prop}

\begin{proof}
A de Bruijn cycle of order $d$ is a binary sequence $\{S(n)\}_{n=1}^{q^d}$ in which every binary $d$-word
appears in a ``window'' $(S(n+1),\ldots,S(n+d))$ for some $j$ (indices taken modulo $q^d$). A \textit{reduced}
de Bruijn cycle is a string of length $q^d-1$ which achieves every $d$-word in some window, except for the word
$0^d$. Note that a reduced de Bruijn cycle may always be turned into an ordinary de Bruijn cycle by inserting an
extra ``0'' into its longest run of 0's.

If $\cP$ is primitive, then $q$ is a generator of $\field{2^d}^\times$, and it is a classical result that the
sequence of constant coefficients of the powers of a multiplicative generator yield a reduced de Bruijn cycle.
Thus, by Proposition~\ref{lem:lsr} the first $2^d-1$ coefficients of $\bar\cP$ are a reduced binary de Bruijn
cycle of order $d$. The reader wishing to explore de Bruijn cycles further can find the basics
in~\cites{MR652466,ShiftRegisterBook}.

Since every string except $0^d$ appears in $\bar \cP$, there are exactly $2^{d-1}$ ones in any period.
\end{proof}

\subsection{Polynomials with non-high density reciprocals}\label{sec:companions}

We see in Figure~\ref{PolyInversePic} that polynomials typically have reciprocals with density near $1/2$. In
Figure~\ref{fig:deltaF,F*}, it is apparent that there is a connection between the density of $\bar\cP$ and
$\bar\cP^\ast$. Our next theorem elucidates the connection.

\begin{figure}[ht]
\begin{center}
\begin{picture}(200,200)
    \put(10,11){\includegraphics[width=180pt]{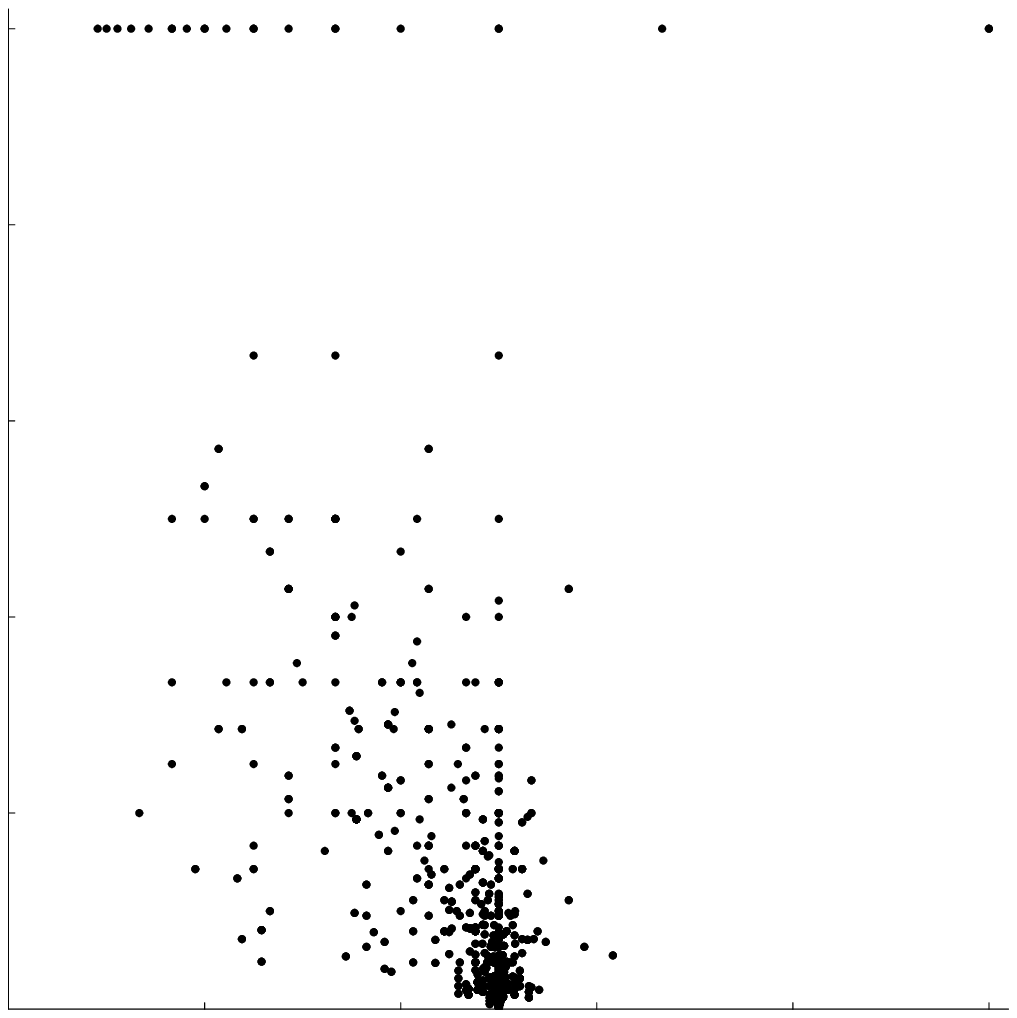}}
    \put(0,0){0}
    \put(90,0){$1/2$}
    \put(183,0){1}
    \put(0,96){$\tfrac12$}
    \put(0,183){1}
\end{picture}
\caption{The points $(\delta(\bar \cP_n),\delta(\bar\cP_n^\ast))$ for odd $n<2^{12}$\label{fig:deltaF,F*}}
\end{center}
\end{figure}

\begin{prop} \label{prop:1/2} If $\ord(\cP)\ge4$, then
    $
    \min\{\delta(\bar\cP),\delta(\bar \cP^\ast)\} \leq 1/2.
    $
\end{prop}

This proposition is best possible in that $\cP_{51}\cP_{15}=1+q^8$, and
$\delta(\bar\cP_{51})=\delta(\bar\cP_{15})=1/2$.

\begin{proof}
Set $D:=\ord(\cP)$. We assume without loss of generality that $\deg(\cP)\le D/2 \le \deg(\cP^\ast)$. If
$\deg(\cP)<3$, then we appeal to the following table of calculations:
    \[{\renewcommand{\arraystretch}{1.5}
    \begin{array}{|c|c|c|c|}\hline
    \cP     &   \cP^\ast                &   \bar \cP^\ast                           &   \delta(\bar \cP^\ast)   \\
    \hline\hline
    1       &   1+q^D                   &   \sum_{n=0}^\infty q^{nD}                &   1/D                     \\
    1+q     &   \sum_{n=0}^{D-1} q^n    &   \sum_{n=0}^\infty (q^{nD}+q^{nD+1})     &   2/D                     \\
    1+q^2   &   \sum_{n=0}^{D-1} q^{2n} &   \sum_{n=0}^\infty (q^{2nD}+q^{2nD+2})   &   1/D                     \\
    1+q+q^2 &(1+q)\sum_{n=0}^{D/3-1}q^{3n}&\sum_{n=0}^\infty (q^{nD}+q^{nD+1}+q^{nD+2})& 3/D        \\ \hline
    \end{array}
    }\]
In the case $\cP=1+q+q^2=\frac{1+q^3}{1+q}$, we see also that $D \equiv 0 \pmod 3$, and by hypothesis $D\ge 4$,
so that $D/3 \le 1/2$.

We assume now that $\deg(\cP)\ge 3$. Since
    \[
    \bar \cP^\ast = \frac{\cP}{1+q^D} = \cP+q^D\cP + q^{2D}\cP + q^{3D} \cP + \cdots
    \]
and $\deg(\cP)<D$, we have $\delta(\bar\cP^\ast) = \ell(\cP)/D.$ If $\cP$ has any zero coefficients, then
$\ell(\cP)\le \deg(\cP)\le D/2$ and so $\delta(\bar\cP^\ast) \leq 1/2$. If $\cP$ has no zero coefficients, then
$\cP = (1+q^{\deg(\cP)+1})/(1+q)$, in which case $\bar \cP = (1+q)/(1+q^{\deg(\cP)+1})$, a series which has
density $2/(\deg(\cP)+1)$. Since $\deg(\cP)\ge 3$, this quantity is $\leq 1/2$.
\end{proof}

\begin{cor}
If $\cP \not \in \{1, 1+q, 1+q+q^2\}$ is a polynomial and $\cP^\ast$ is primitive, then $\delta(\bar \cP) \leq
1/2$.
\end{cor}

\subsection{Eventually periodic sets}\label{sec:eps}

An eventually periodic set is one whose generating function has the form $\cE(q)+\frac{\cP(q)}{1+q^D}$, for some
polynomials $\cE, \cP$ with $\deg(\cP)<D$, and exactly one of $\cE,\cP$ has constant term 1. The finite sets
containing 0 are examples. Another example is given by the set $\NN\setminus\{n\colon n\equiv 2 \pmod{4}\}$
(which has density $3/4$), whose reciprocal is the set $\{1\} \cup \{n \in \NN \colon n \text{ is congruent to
0, 2, 5, or 6 modulo 7}\}$ (which has density $4/7$).

\begin{prop}
The reciprocal of an eventually periodic set is an eventually periodic set.
\end{prop}

This proposition is essentially the same as that which asserts that rational numbers have eventually periodic
decimal expansions.

\begin{proof}
Obviously, the reciprocal of a ratio of polynomials (each with constant term 1) is a ratio of polynomials (each
with constant term 1). All that we need to observe is that such a ratio $\cR/\cS$ can be written in the form
    \[
    \frac{\cR}{\cS} = \cE+\frac{\cQ}{1+q^D},
    \]
with $\deg(\cQ)<D$. By long division, we can write $\cR/\cS$ in the form $\cE+\cP/\cS$ with
$\deg(\cP)<\deg(\cS)$. But this is the same as $\cE+\frac{\cP \cS^\ast}{1+q^D}$, where $\cS \cS^\ast = 1+q^D$,
and $\deg(\cP \cS^\ast)<\deg(\cS \cS^\ast)=D$.
\end{proof}

\section{The powers of two}\label{sec:PowersOf2}

We saw in Section \ref{sec:Polynom} that the reciprocal of a polynomial (other than $\cP_1$) has positive
density. One might wonder if the reciprocal of any set with zero density has positive density. Our next theorem
shows that this is not the case.

We note the $m=1$ case of Theorem~\ref{thm:PowersOf2}: the reciprocal of $A_1=\{0\}\cup \{2^n \colon n\in\NN\}$
is $\bar A_1=\{2^n-1\colon n\in\NN\}$. This is easily proved directly by considering the following
sum-preserving involution on $A_1\times \bar A_1$. For $s,t\in\NN$ and distinct, define
    $\mu(0,0)=(0,0)$,
    $\mu(0,2^{t+1}-1) = (2^t, 2^t-1)$,
    $\mu(2^s,2^t-1) = (2^t,2^s-1)$,
    $\mu(2^t,2^t-1) = (0, 2^{t+1}-1)$.
The existence of this sum-preserving fixed-point-free involution proves that every positive integer $n$ can be
written in the form $a+\bar a$, where $a\in A_1$ and $\bar a\in\bar A_1$, in an even number of ways. A similar
proof can be given for $m=2$, and presumably for any $m$, but quickly grows tedious. We now give an algebraic
proof that does not depend on $m$.

\begin{thm}\label{thm:PowersOf2}
Let $m\ge1$. The reciprocal of the set $A_m:=\{0\}\cup\{2^{mn} \colon n\in \NN\}$ is the set
    \[
    \bar A_m:= \bigg\{-1+\sum_{i=0}^{m-1} x_i 2^{i+mn_i} \colon x_i\in\{0,1\}, \vec{x}\not=\vec{0}, n_i \in \NN \bigg\}.
    \]
In particular, both $\delta(A_m,n)$ and $\delta(\bar A_m,n)$ are $O_m\big({\frac{\log n}n}\big)$.
\end{thm}

\begin{proof}
Set $\cF(q)=\sum_{n\ge0} q^{2^{mn}}$. By the children's binomial theorem\, $\cF(q^2)=\cF(q)^2$, and consequently
by induction we see that $\cF(q^{2^m})=\cF(q)^{2^m}$.

Now, by the definition of $\cF$, $\cF(q^{2^m})=\cF(q)+q$ and so
    \begin{align*}
    q   &=\cF(q)^{2^m}+\cF(q)\\
        &=\big(1+\cF(q)\big)\, \left(\cF(q)+\cF(q)^2+\cF(q)^3+\dots + \cF(q)^{2^m-1} \right)\\
        &=\big(1+\cF(q)\big)\,\left(1+\prod_{i=0}^{m-1} \big(1+\cF(q)^{2^i}\big) \right) \\
        &=\big(1+\cF(q)\big)\,\left(1+\prod_{i=0}^{m-1} \big(1+\cF(q^{2^i})\big) \right)
    \end{align*}
The series $1+\cF(q)$ is the generating function of $\{0\}\cup \{2^{mn}\colon n\in\NN\}$, and
\mbox{$1+\prod_{i=0}^{m-1} (1+\cF(q^{2^i}))$}\ is the generating function of \mbox{$\big\{\sum_{i=0}^{m-1} x_i
2^{i+mn_i} \colon x_i\in\{0,1\}, \vec{x}\not=\vec{0},n_i\in\NN\big\}$}, so this identity is equivalent to the
theorem.
\end{proof}

The reader may be interested to note that the reciprocal of the extremely thick set $\NN\setminus\{2^n\colon
n\in\NN\}$ is the thin set $\{0,3\}\cup\{2^n-1,2^n-3\colon n\ge3\}$, whereas the reciprocal of
$\NN\setminus\{4^n\colon n\in\NN\}$ appears to have density $1/2$.

Our next theorem shows that the examples given by Theorem~\ref{thm:PowersOf2} are extremal. It is impossible for
a set and its reciprocal to both grow sub-logarithmically. This result was suggested to us by Ernest Croot
[personal communication].

\begin{thm}
Let $F,\bar F$ be reciprocals (not both $\{0\}$), and suppose that $r$ is the least positive integer in
$F\cup\bar F$. Then
    \[\big| F \cap [0,n]\big| + \big| \bar F \cap [0,n] \big| \ge 2+\floor{\log_2 (n/r)}.\]
\end{thm}

\begin{proof}
First, note that $r\in F\cap\bar F$. Let $N\ge r$, so that neither $F\cap[1,N)$ nor $\bar F \cap[1,N)$ is empty,
and let $m,\bar m$ be the largest elements of those sets. Since $q^{m+\bar m}$ occurs in the product $\cF \bar
\cF$ at least once, it must occur at least twice. Since $N\le m+\bar m < 2N$, we see that
    \[\big| F \cap [N,2N)\big| + \big| \bar F \cap [N,2N) \big| \ge 1.\]
Straightforward counting concludes the proof, since $F\cup\bar F$ contains 0 twice, and must intersect each of
the intervals $[r,2r)$, $[2r,2^2r)$, $[2^2r,2^3r)$, $\dots$.
\end{proof}

\section{Theta functions}\label{sec:Theta}

Every quadratic that takes integers to integers can be written in the form
    $
    c_0+c_1 n + c_2 \frac{n(n-1)}2
    $
with $c_i\in\ZZ$. We wish to study the ranges of such quadratics, but we only wish to consider sets that contain
0; without loss of generality we may take $c_0=0$. Thus, we set
    \[
    \Theta(c_1,c_2) := \left\{ c_1 n+ c_2 \frac{n(n-1)}2  \colon n\in \ZZ \right\}.
    \]
Moreover, we are only interested in those sets that consist of nonnegative integers, so we may assume that
$c_2\ge c_1 \ge 0$. And since $\Theta(c_1,c_2)=\Theta(c_2-c_1,c_2)$ we may also assume that $c_2\ge 2c_1$.
Finally, we are only interested in those sets whose $\gcd$ is 1: we can assume that $\gcd(c_1,c_2)=1$. The only
set with $c_1=0$ not excluded is $\Theta(0,1)=\{\tbinom n2 \colon n\ge 1\}$, and the only set with $c_2=2c_1$
that is not excluded is $\Theta(1,2)=\{n^2 \colon n\ge 0\}$. Otherwise, we have $c_2>2c_1>0$, and
$\gcd(c_1,c_2)=1$.

In Figure~\ref{theta.data.table}, we give the number of elements in the reciprocal of $\Theta(c_1,c_2)$ (with
$c_2\le 18$) that are at most $10^5$. We note that none of the entries of this table are larger than 50450, and
the entries that are less than 49750 are exactly those with $c_2\equiv 2\pmod{4}$. This computation partially
justifies Conjecture~\ref{c1c2Conjecture}.

\begin{figure}[t]
\begin{center}
$c_1$ \vskip4pt $c_2$ \quad
\begin{tabular}{c||ccccccccc}
     &0 &  1    &  2    &  3    &  4    &  5    &  6 & 7 & 8   \\ \hline\hline
 1  &50162 &   &   &   &   &   & \\
 2  &   &17317 &   &   &   &   &   \\
 3  &   &50201 &   &   &   &   &   \\
 4  &   &50162 &   &   &   &   &   \\
 5  &   &50265 & 49994 &   &   &   &   \\
 6  &   &17814 &   &   &   &   &   \\
 7  &   &50062 & 50187 & 50449 &   &   &   \\
 8  &   &50042 &   & 49944 &   &   &  \\
 9  &   &50214 & 49827 &   & 50023 &   &   \\
 10 &   &34009 &   & 36084 &   &   &   \\
 11 &   &49918 & 50181 & 49918 & 49943 & 49856 &   \\
 12 &   &49869 &   &   &   & 50254 &   \\
 13 &   &50089 & 49752 & 49988 & 49992 & 50295 & 49912 \\
 14 &   &40981 &   & 41776 &   & 39062 &   \\
 15 &   &50004 & 50195 &   & 49949 &   &  & 49900 \\
 16 &   &50001 &   & 49924 &   & 49996 &  & 50090 \\
 17 &   &50198 & 49921 & 49932 & 50052 & 50114 & 49826 & 49818 & 49816\\
 18 &   &48224 &   &   &   & 44500 &  & 43772 \\
\end{tabular}
\end{center}
\caption{The number of elements $\le 100000$ in the reciprocal of $\Theta(c_1,c_2)$\label{theta.data.table}}
\end{figure}

There is another property of $\Theta(c_1,c_2)$ that happens exactly when $c_2\equiv 2\pmod{4}$: the set
$\Theta(c_1,c_2)$ is not uniformly distributed modulo 4.

\begin{prop}
Let $\gcd(c_1,c_2)=1$. The set $\Theta(c_1,c_2)$ is uniformly distributed modulo every power of 2 if and only if
$c_2 \not\equiv 2\pmod{4}$.
\end{prop}

\begin{proof}
First, suppose that $c_1=2k+1$ and $c_2=4\ell+2$. Set
    \[
    f(n):= c_1 n + c_2\frac{n(n-1)}{2} = (2\ell+1)n^2+2(k-\ell)n.
    \]
If $k$ and $\ell$ have the same parity, then $f(n)\equiv (2\ell+1) n^2 \pmod{4}$, and since $n^2$ takes only two
values modulo 4, the set is not uniformly distributed modulo 4. If $k$ and $\ell$ have different parity, then
    \[
    (2\ell+1)n^2+2(k-\ell)n \equiv (2\ell+1)n^2+2n \pmod{4}
    \]
only takes on the values $0,3$ modulo 4. Thus, if $c_2\equiv 2 \pmod 4$, then $\Theta(c_1,c_2)$ is not uniformly
distributed modulo 4.

Now suppose that $c_2=4\ell$, and since $\gcd(c_1,c_2)=1$, we know that $c_1$ is odd. We have
    \[
    f(n):=c_1 n + c_2\frac{n(n-1)}{2} =  2\ell n^2+(c_1-2\ell)n\equiv n \pmod{2}.
    \]
The formal derivative of $f(n)$ is $4\ell n + c_1-2\ell \not\equiv 0 \pmod{2}$. By Hensel's
Lemma\footnote{Hensel's Lemma: If $f(n)$ is a polynomial with integer coefficients, and the two congruences
$f(n)\equiv a \pmod{p}, f'(n)\not\equiv 0 \pmod{p}$ have a simultaneous solution, then $f(n)\equiv a$ has a
unique solution modulo every power of the prime $p$.}, the range of the polynomial $f(n)$ hits every congruence
class modulo every power of $2$. Since for every $j$, $f(n)$ is periodic modulo $2^j$ with period $2^j$, we see
that it is uniformly distributed modulo $2^j$.

Now suppose that $c_2=2\ell+1$ is odd. Set
    \begin{align*}
    G &:= \left\{ (2\ell+1)m(2m-1)+c_1(2m) \colon m\in \ZZ\right\} \\
    H &:= \left\{ (2\ell+1)(2m+1)m+c_1(2m+1) \colon m\in \ZZ\right\}
    \end{align*}
so that $\Theta(c_1,2\ell+1)=G\cup H$. The set $G$ is the range of $g(m):=f(2m)=(2\ell+1)m(2m-1)+2c_1m \equiv m
\pmod{2}$, which has derivative $g'(m)\equiv 1 \pmod{2}$, and the set $H$ is the range of
$h(m):=(2\ell+1)(2m+1)m+c_1(2m+1)\equiv m+c_1 \pmod{2}$, which has derivative $h'(m)\equiv 1 \pmod{2}$. Thus, by
Hensel's Lemma, both $G$ and $H$ exhaust every congruence class modulo $2^j$, and by periodicity of $g(m)$ and
$h(m)$ are therefore uniformly distributed modulo $2^j$.
\end{proof}

\subsection{The squares}\label{sec:squares}

Let ${\cS}(q)=\sum_{n=0}^\infty q^{n^2}$, and $S=\{0,1,4,9,16,25, \dots\}$. Figure~\ref{SquarePics} shows
$\delta(\bar S,x)$ for two ranges of $x$. On the small scale, we see that the relative density behaves
irregularly, with many small increases and decreases. On the larger scale, we see that the relative density
seems to decrease inexorably.

We characterize completely the values of $\bar S$ in the residue classes $0, 1, 2 \pmod{4}$.

Let $\nu_p(n)$ be the integer such that $p^{\nu_p(n)} \mid n$ and $p^{\nu_p(n)+1} \nmid n$, so that
    \[n=\prod_{\text{$p$ prime}} p^{\nu_p(n)}\]
for every $n$. Let $r_2(n)$ be the number of representations of $n$ in the form $y^2+z^2$, where $y$ and $z$ are
integers.

\begin{thm}\label{thm:SquareDescription} Let $n\in\NN$. If $n$ is even, then $n\in\bar S$ if and only if $n$ is
twice a square. If $n\equiv1\pmod4$ is not a square, then $n\in\bar S$ if and only if $\nu_p(n)$ is even for
every prime $p$ except one, and that prime $p$  and $\nu_p(n)$ are both congruent to 1 modulo 4. If $n\equiv
1\pmod 4$ is a square, then $n\in\bar S$ if and only if $\nu_p(n)\equiv 2 \pmod4$ for an even number of primes
$p\equiv 1 \pmod 4$.
\end{thm}

We will need the following lemmas. The first expresses $\bar s_n$ in terms of the number of representations of
$n$ by a particular (depending on $n$) quadratic form. The second is quoted without proof from~\cite{MR1083765},
and gives a formula for $r_2(n)$.

\begin{lem}\label{lem:quadraticform}
Let $n\in\NN$, and let $j\in\NN$ satisfy $n\equiv 2^j-1 \pmod{2^{j+1}}$. Then $\bar s_n=1$ if and only if
    \[
    \#\bigg\{(k_0, \ldots,k_{j-1}, k_{j+1}) \colon \,k_i\in\NN,\,
        n= 2^{j+1} k_{j+1}^2+\sum_{i=0}^{j-1} 2^i k_i^2\bigg\}
    \]
is odd.
\end{lem}

\begin{proof}
By Lemma~\citebarfnitem{ii}, $\bar s_n=1$ exactly if there are an odd number of tuples $(k_0,k_1,\dots)$ with
weight
    \begin{equation}\label{equ:representation}
    n=k_0^2+2 k_1^2+4 k_2^2+8 k_3^2+\cdots.
    \end{equation}
Let $w(n)$ be the number of such tuples. We give a weight-preserving involution $\mu$ of such tuples, and deduce
the lemma from
    \[
    w(n)\equiv \#(\text{fixed points of $\mu$ with weight $n$}) \pmod{2}.
    \]

Since $n\not\equiv 2^i-1\pmod{2^{i+1}}$ for $0\le i < j$, reducing~\eqref{equ:representation} modulo
$2,4,\ldots,2^j$ successively tells us that $k_0, k_1,\ldots, k_{j-1}$ are odd, while $n\equiv 2^j-1
\pmod{2^{j+1}}$ tells us that $k_j$ is even. Now define $J$ to be the least integer with the two properties:
$J\ge j+2$; and $2k_J\not=k_j$.

We define
    \[
    \mu(k_0,k_1,k_2,\dots) = (k_0,k_1,\ldots, k_{j-1},2k_J,k_{j+1},k_J,k_J,\ldots,k_J,k_j/2,k_{J+1},k_{J+2},\ldots),
    \]
where $k_J$ is repeated $J-j-2$ times. That this is a weight-preserving involution is a routine calculation.

The fixed points of $\mu$ are those tuples with $0=k_j=k_{j+2}=k_{j+3}=\cdots$. In other words, there is a fixed
point for each solution to
    \[
    n=k_0^2+2k_1^2+\cdots+2^{j-1}k_{j-1}^2+2^{j+1}k_{j+1}^2.\hfill \qed
    \]
\renewcommand{\qed}{}
\end{proof}

\begin{lem}[\cite{MR1083765}*{Theorem 3.22}]\label{lem:fromNMZ}
If $\nu_p(n)$ is odd for any prime $p$ congruent to 3 (modulo 4), then $r_2(n)=0$. Otherwise, $r_2(n)= 4
\prod_{p} (\nu_p(n)+1)$, where the product extends over all primes congruent to 1 (modulo 4).
\end{lem}

\begin{proof}[Proof of Theorem~\ref{thm:SquareDescription}]
If $n$ is even, then $n\equiv 2^0-1 \pmod{2^{0+1}}$, so we can apply Lemma~\ref{lem:quadraticform} with $j=0$ to
arrive at $\bar s_n=1$ if and only if $n$ has an odd number of representations of the form $2k_1^2$ (with
$k_1\ge 0$). Clearly there cannot be more than one such representation, and there is one exactly if $n$ is twice
a perfect square.

If $n\equiv 1\pmod 4$, then we may apply Lemma~\ref{lem:quadraticform} with $j=1$ to arrive at $\bar s_n=1$ if
and only if $n$ has an odd number of representations of the form $k_0^2+4k_2^2$ (with $k_0$ and $k_2$
nonnegative).

We assume for now that $n$ is not a square. Since $n$ is odd, there are no such representations with $k_0=0$,
and since $n$ is not a square, there are no such representations with $k_2=0$. Thus, every such representation
$k_0^2+4 k_2^2$ gives rise to 8 representations $\{(\pm k_0)^2+(\pm 2 k_2)^2, (\pm 2k_2)^2+(\pm k_0)^2\}$ of $n$
in the form $y^2+z^2$. Moreover, any solution to $n=y^2+z^2$ must have one of $y$ or $z$ even and the other odd
since $n$ is odd, and $y\not= z$ since $n$ is odd. Since $n$ is not a square, neither $y$ nor $z$ is zero. Every
representation $(y,z)$ occurs as one of a family of 8 such representations, and one of these has
$n=y^2+z^2=y^2+4 (z/2)^2$ with $y>0$ and $z>0$. Thus, $\bar s_n=1$ if and only if $r_2(n)/8$ is odd.

By Lemma~\ref{lem:fromNMZ}, $r_2(n)/8=0$ if $\nu_p(n)$ is odd for any prime $p$ congruent to 3 modulo 4.
Otherwise, $r_2(n)/8=\frac 12 \prod_p (\nu_p(n)+1)$, where the product extends over those primes that are
congruent to 1 modulo 4 (in the remainder of this paragraph, $p$ is always 1 modulo 4). First, note that
$\nu_p(n)$ is odd for some prime $p$ since $n$ is not a square. If some $\nu_p(n)$ is 3 modulo 4 for some $p$,
then $r_2(n)/8$ is even, and similarly if $\nu_p(n)$ is 1 modulo 4 for two primes $p$. Thus, $r_2(n)/8$ is odd
precisely if $\nu_p(n)$ is odd for exactly one prime, and both that prime and $\nu_p(n)$ are 1 modulo 4.

Now we assume that $n\equiv 1\pmod4$ is a square, say $n=x^2$. Then, as above, most representations of $n$ in
the form $k_0^2 + 4 k_2^2$ correspond to 8 representations of $n$ in the form $y^2+z^2$, but the representation
$n=x^2+4\cdot 0^2$ only corresponds to 4 representations in the form $y^2+z^2$. Since $n$ is a square, we know
that $\nu_p(n)$ is even for every prime $p$. Thus, $\bar s_n=1$ if and only if
    \[\frac{r_2(n)-4}8+1\]
is odd. Using the formula from Lemma~\ref{lem:fromNMZ}, this happens exactly if $1\equiv \prod_p
(\nu_p(n)+1)\pmod4$, where the product extends over primes that are 1 modulo 4. This, in turn, happens exactly
when $\nu_p(n)\equiv 2 \pmod4$ for an even number of primes $p\equiv1\pmod4$.
\end{proof}

We suspect that $\delta(\bar S)=0$ and that this may follow from the theory of modular forms, but again, this is
outside the scope of this paper. We emphasize in Corollary~\ref{cor:barSdensity} that our characterization of
$\bar S$ is consistent with Conjecture~\ref{c1c2Conjecture}.

\begin{cor}\label{cor:barSdensity} The set $\{n\in\NN \colon n\in \bar S, n\not\equiv 3 \pmod4\}$ has zero
density.
\end{cor}

\begin{proof}
By Theorem~\ref{thm:SquareDescription}, the set $\bar S$ clearly has no density in $0\bmod 2$. We will use the
description given in Theorem~\ref{thm:SquareDescription} to show that $\bar S$ also has zero density in $1\bmod
4$.

By the Wiener-Ikehara Theorem (see~\cite{MR2111739}*{Section 7.2}), we have for any set $A$ of positive integers
    \[
    \lim_{n\to\infty} \delta(A,n) = \lim_{s\to 1+} (s-1)\sum_{a\in A} a^{-s}.
    \]
Set $A=\{n^2 p \colon 1\le n\in \NN, p \text{ prime}\}$, and observe that $\delta(A)=0$ since
    \begin{align*}
    \delta(A) &\leq \lim_{s \to 1+} (s-1) \sum_{a\in A} a^{-s}\\
    &= \lim_{s \to 1+} (s-1) \left(\prod_{p \text{ prime}} (1-p^{-2s})^{-1}\right)\left(\sum_{p\text{ prime}}
    p^{-s}\right)\\
    &= \lim_{s \to 1+} (s-1) \zeta(2s) \left(\sum_{p\text{ prime}}     p^{-s}\right)\\
    &= \zeta(2) \;\delta(\text{primes}) = 0
    \end{align*}

Note that the subset of $\bar S$ whose elements are even has density 0, and the subset whose elements are
congruent to 1 modulo 4 is (except for some squares) contained in $A$. Thus
    \[\delta(\{n\in\bar S\colon n\not\equiv 3\pmod{4}\}) \le\delta(\text{squares})+ \delta(A)=0.\hfill\qed\]
\renewcommand{\qed}{}
\end{proof}

\section{Prouhet-Thue-Morse numbers}\label{sec:Thue}

Set $t_n=1$ if the binary expansion of $n$ contains an even number of ``\texttt{1}''s, and set $t_n=0$
otherwise. The set $T:=\{n \colon t_n=1\}=\{0, 3, 5, 6, 9, \dots\}$ is called the Prouhet-Thue-Morse sequence.
This sequence frequently arises because it simultaneously has enough structure to analyze, and enough
``random-like'' behavior to be interesting. The survey~\cite{MR1843077} details four of the occasions that the
sequence has been independently rediscovered: first in number theory (Prouhet), then combinatorics (Thue), then
in differential geometry (Morse), and finally chess grandmaster Max Euwe rediscovered it to demonstrate that the
rules then in use did not imply that chess is a finite game.

For every $n\in\NN$, $2n\in T$ if and only if $2n+1\not\in T$; thus $\cT(q):=\sum_{n=0}^\infty t_n q^n$ has
$\delta(\cT)=1/2$. The sequence $t_0, t_1, \dots$ is not eventually periodic (in fact, the real number with
binary expansion $0.t_0t_1t_2\cdots$ is transcendental~\cites{MR0457363,MR1869317}), so $\bar \cT$ is not a
polynomial. A counting argument~\cite{MR1843077} reveals the interesting identity:
    \begin{equation}\label{equ:ThueFunctional}
    (1+q)^3\cT(q)^2+(1+q)^2 \cT(q)=q.
    \end{equation}
Multiplying by $\bar \cT(q)$ yields $q\bar \cT(q)=(1+q+q^2+q^3)\cT(q)+1+q^2$, whence for $n\ge2$
    \[
    \bar t_n = t_{n+1}+t_n+t_{n-1}+t_{n-2}.
    \]
This leads reasonably directly (albeit with the modest labor involved in deriving~\eqref{equ:ThueFunctional}) to
a proof of Theorem~\ref{thm:ThueMorse}. Instead, we give a proof which does not rely on the special form of the
functional equation \eqref{equ:ThueFunctional}, and so is more representative of the process of {\em finding}
reciprocals.

\begin{thm}\label{thm:ThueMorse}
The reciprocal of the set $T$ of Prouhet-Thue-Morse numbers is
    \[
    \bar T=\{0\}\cup \{4k\pm1 \colon \text{the binary expansion of $k\ge1$
    ends in an even number of ``\texttt{1}''s}\}.
    \]
    Consequently, $\delta(\bar \cT)=1/3$.
\end{thm}

If (the binary expansion of) $k$ ends in an even number of ``\texttt{1}''s, then $4k+1$ ends with a string
\texttt{10${}^{2k+1}$1} (a ``\texttt{1}'' followed by an odd number of ``\texttt{0}''s followed by a single
``\texttt{1}''), while $4k-1$ ends with a string \texttt{01${}^{2k}$} (a ``\texttt{0}'' followed by a positive
even number of ''\texttt{1}''s).

\begin{proof}
By Lemma~\ref{lem:golden},
$\bar \cT(q)= \sum_{n=0}^\infty r(n) q^n$,
where $r(n)$ is the number of ways to
write $n$ as
    \[
        n=s_0+2s_1+4s_2 + 8 s_3 + \cdots + 2^k s_k + \cdots
    \]
where the $s_k$ are Prouhet-Thue-Morse numbers. We will build an involution $\tau$ on the set of such
representations, and $r(n)$ will have the same parity as the number of fixed points of $\tau$.

By a {\em tuple}, we mean an infinite list of Prouhet-Thue-Morse numbers which is 0 from some point on. The {\em
weight} of a tuple $(s_0,s_1,\ldots)$ is $\sum_{n=0}^\infty s_n 2^n$.

We now give the weight-preserving permutation $\tau$ of the set of tuples which is actually an involution. The
permutation $\tau$ has an odd number of fixed points with weight $n>0$ if and only if the binary expansion of
$n$ ends with a string \texttt{10${}^{2k+1}$1} (a ``\texttt{1}'' followed by an odd number of ``\texttt{0}''s
followed by a single ``\texttt{1}'') or ends with a string \texttt{01${}^{2k}$} (a ``\texttt{0}'' followed by a
positive even number of ''\texttt{1}''s). These are exactly the numbers of the form $4k\pm1$, where the binary
expansion of $k$ ends in an even number of ``\texttt{0}''s, and this will conclude the proof.

\paragraph{Defining the permutation $\tau$:}

Suppose that $s_0$ is even. If $s_0\not=2s_1$, then set
    \[
    \tau(s_0,s_1,s_2,\ldots) := (2s_1,s_0/2,s_2,s_3,\dots).
    \]
If $s_0=2s_1$, then let $i$ be minimal with $s_1\not=s_i$, and set
    \[
    \tau(s_0,s_1,s_2,\ldots) :=
    (2s_i, s_i, s_i, \ldots, s_i,s_1,s_{i+1},s_{i+2},\ldots),
    \]
where $s_i$ is repeated $i-1$ times.
The only fixed point with $s_0$ even is $(0,0,\ldots)$ with weight 0.

Now suppose that $s_0\equiv 3 \pmod{4}$.
Since $(s_0-3)/2$ is even, we can define $v_0,v_2,v_3,\dots$ by
    \[
    (v_0,v_2,v_3,\dots):= \tau((s_0-3)/2, s_2,s_3, \ldots),
    \]
where the action of $\tau$ has already been defined above. Note that $v_1$ is not defined, and that $s_1$ has
not been used. We now set
    \[
    \tau(s_0,s_1,s_2,\ldots) := (2v_0+3,s_1,v_2,v_3, \dots).
    \]
The only fixed points with $s_0\equiv 3 \pmod 4$ are the tuples of the form $(3,s_1,0,0,\dots)$, where $s_1$ is
a Prouhet-Thue-Morse number. These fixed points have weight $3+2s_1$.

Now suppose that $s_0 \equiv 1 \pmod 4$. If there exists an $L$ such that $s_i=(s_0+1)/2$ for $0<i\le L$ and
$s_i=0$ for $i> L$, then we let $\tau$ fix the tuple. These will be the only fixed points of $\tau$ with
$s_0\equiv 1 \pmod 4$, and will have weight $2^L s_0+2^L-1$. Otherwise, if any $s_i$ is even (except for the
tail of zeros in the tuple $(s_0,s_1,s_2,\ldots)$), then let $K:=\min\{i \colon s_i \text{ even}\}$, and set
    \[
    \tau(s_0,s_1,s_2,\ldots)
    := (s_0, s_1, s_2, \ldots, s_{K-1}, \tau(s_K, s_{K+1}, s_{K+2},\dots)).
    \]
If on the other hand all $s_i$ are odd (except for the ending string of zeros),
then define $v_0, v_1, v_2,\ldots$ by
    \[
    (v_0, v_1, v_2, \ldots):= \tau(s_0+1, s_1, s_2, s_3, \ldots)
    \]
and set
    \[
    \tau(s_0,s_1,s_2,\ldots) := (v_0-1, v_1, v_2, \ldots).
    \]

That $\tau$ is an involution with precisely the claimed fixed points
is simply a matter of checking the various cases;
we cheerfully leave this important tedium to the reader.

\paragraph{Analysis of $\tau$'s fixed points with weight $n$:}

Suppose that $n$ is even. By parity considerations, we see that all tuples $(s_0, s_1, \ldots)$ with weight
$\sum_{i=0}^\infty s_i 2^i=n$ have $s_0$ even. Since the only fixed point with $s_0$ even is $(0,0,\ldots)$, we
see that $r(0)=1$ and $r(n)$ is even for all even $n>0$. From this point on we assume that $n$ is odd.

Suppose that $n\equiv 1\pmod 4$, and $(s_0,s_1, \ldots)$ is a fixed point of $\tau$ with weight $n$. Since $n$
is odd, $s_0$ is either 1 or 3 modulo 4. If $s_0\equiv 1 \pmod{4}$, then $s_1\equiv 1 \pmod{2}$, and such a
tuple can be fixed by $\tau$ only if $n=s_0$, and $n$ is a Prouhet-Thue-Morse number. If $s_0\equiv 3 \pmod 4$,
then $s_1\equiv 1 \pmod{2}$, and such a tuple can be fixed by $\tau$ only if $n=3+2s_1$, i.e., if $(n-3)/2$ is a
Prouhet-Thue-Morse number (and in this case there is exactly one such tuple). Thus $\tau$ has either 0, 1, or 2
fixed points, and we care about when it has an odd number of fixed points. Since $n\equiv 1 \pmod 4$, the binary
expansion of $n$ can be written as $(\texttt{x10${}^k$1})_2$ for some binary string \texttt{x} and positive
integer $k$. We see that the binary expansion of $(n-3)/2$ is $(\texttt{x01${}^k$})_2$. Thus, if $k$ is even,
then either both $n$ and $(n-3)/2$ are Prouhet-Thue-Morse numbers or neither is. If $k$ is odd, then exactly one
of $n$ and $(n-3)/2$ are Prouhet-Thue-Morse numbers. Hence, $\tau$ has an odd number of fixed points exactly if
the binary expansion of $n$ ends in $\texttt{10${}^k$1}$, with $k$ an odd number.

Now suppose that $n\equiv 3 \pmod 4$, and $(s_0,s_1, \ldots)$ is a fixed point of $\tau$ with weight $n$. Since
$n$ is odd, $s_0$ is either 1 or 3 modulo 4. If $s_0\equiv 1 \pmod{4}$, then $s_1\equiv 1 \pmod{2}$, and such a
tuple can be fixed by $\tau$ only if $s_i=(s_0+1)/2$ for all $0<i\le L$ and $s_i=0$ for $i>L$. In this case,
$n=2^L s_0+(2^L-1)$. Since $s_0\equiv1 \pmod 4$, this implies that the binary expansion of $n$ ends with $L+1$
``\texttt{1}''s (in particular, at most one value of $L$ can lead to such a fixed point). Moreover, $2^L s_0 +
(2^L-1)$ is a Prouhet-Thue-Morse number if and only if $L$ is even. If $s_0\equiv 3 \pmod4$, then $s_1\equiv 0
\pmod 2$, and such a tuple is fixed if and only if it is of the form $(3, s_1, 0,0,\ldots)$. This can happen
exactly if $(n-3)/2$ is a Prouhet-Thue-Morse number.

Suppose that $n$ is a Prouhet-Thue-Morse number. If the binary expansion of $n$ ends in exactly $2k>0$ ``1''s,
then $(3,(n-3)/2,0,0,\ldots)$ is the only fixed point of $\tau$. If the binary expansion of $n$ ends in $2k+1>0$
``\texttt{1}''s, then both $(3,(n-3)/2,0,0,\ldots)$ and
    \[
    \left( \frac{n-2^{2k}+1}{2^{2k}}, \frac{n+1}{2^{2k+1}}, \frac{n+1}{2^{2k+1}}, \ldots, 0,0,\ldots \right)
    \]
(the term $(n+1)/2^{2k+1}$ is repeated $2k$ times) are fixed points.

Now suppose that $n$ is not a Prouhet-Thue-Morse number. If the binary expansion of $n$ ends in exactly $2k>0$
``\texttt{1}''s, then
    \[
    \left( \frac{n-2^{2k-1}+1}{2^{2k-1}}, \frac{n+1}{2^{2k}}, \frac{n+1}{2^{2k}}, \ldots, 0,0,\ldots \right)
    \]
(the term $(n+1)/2^{2k}$ is repeated $2k-1$ times) is the only fixed point. If the binary expansion of $n$ ends
in $2k+1>0$ ``\texttt{1}''s, then there are no fixed points.
\end{proof}



\end{document}